\documentclass[red,11pt,a4paper]{article}

\usepackage{cite}

\usepackage{amsmath}
\usepackage{amscd}
\usepackage{amssymb}
\usepackage[pdftex]{color,graphicx,hyperref}
\usepackage{latexsym}
\usepackage{color}
\usepackage[normalem]{ulem}

\bf
\usepackage{ika}

\title{Two-velocity hydrodynamics in fluid mechanics: Part I \\
Well posedness for zero Mach number systems}
\author{Didier Bresch$^1$, Vincent Giovangigli$^{2}$, Ewelina Zatorska$^{2,3,4}$}

\begin{document}
\maketitle
\normalsize
\begin{center}
{\small 1.  Universit\'e de Savoie, Laboratoire de Math\'ematiques\\
73376 Le Bourget du Lac, France\\
      {\ }\\
2.  Centre de Math\'{e}matiques Appliqu\'{e}es,\\
      \'{E}cole Polytechnique, 91128 Palaiseau Cedex, France.\\
   {\ }\\
3. Institute of Mathematics\\
Polish Academy of Sciences, ul \'Sniadeckich 8, 00-656 Warszawa, Poland\\
      {\ }\\
4. Institute of Applied Mathematics and Mechanics\\
University of Warsaw, ul. Banacha 2, 02-097 Warszawa, Poland\\}
\end{center}

\noindent{\bf Abstract:} In this paper we prove  global in time existence of weak solutions to zero Mach number systems arising in fluid mechanics with periodic boundary conditions.
   Relaxing a certain algebraic constraint between the viscosity and the conductivity introduced in  [{\sc D.~Bresch, E.H. Essoufi,} and {\sc M. Sy}, \it J. Math. Fluid Mech. \rm  2007] gives a more complete answer to an open question formulated in [{\sc P.--L.~Lions}, Oxford 1998].
     We introduce a new mathematical entropy which clearly shows existence of two-velocity hydrodynamics with a fixed mixture ratio.
As an application of our result we first discuss a model of gaseous mixture extending the results of [{\sc P. Embid}, {\it Comm. Partial Diff. Eqs}. 1987] to the global weak solutions framework. Second, we present the ghost effect system studied by [{\sc C.D.~Levermore, W.~Sun, K.~Trivisa}, {\it SIAM J. Math. Anal}. 2012] and discuss a contribution of the density-dependent heat-conductivity coefficient to the issue of existence of weak solutions. 

\medskip


\section{Introduction} The low Mach number limit for classical solutions to the full compressible Navier-Stokes  equations was studied notably by {\sc T. Alazard} in  \cite{Al06}. When the large temperature variations and thermal conduction are taken into account, the limit system reads (see (1.7) page 7 in \cite{Al06})
\begin{equation}\label{Lim}
\begin{array}{c}
\gamma P_0  \Div \vu = -(\gamma -1) \kappa \Div\vc{Q},\\
\vr \lr{\pt \vu+ \vu\cdot \nabla \vu}
         +\nabla \pi =- \delta \Div \vc{S} ,\\
\vr C_P \lr{\partial_t T + \vu \cdot\nabla T} =- \kappa \Div\vc{Q},
\end{array}
\end{equation}
where the unknown $\vr,\vu,\pi,T$ denote the fluid density, velocity vector field, pressure and temperature, respectively, $\vc{S}$ and $\vc{Q}$ denote the viscous tensor and the heat flux.\\
The two dimensionless parameters distinguished in  \cite{Al06}
$$ \delta \in [0,1], \qquad \kappa \in [0,1],$$
are the inverse of the {\it Reynolds} number and the {\it P\'eclet} number, respectively, measuring the importance of the viscosity and the heat-conduction.\\
 The density and temperature of the fluid are related by
\eq{\vr = P_0/(RT)\label{press},} 
where $P_0$ denotes the constant pressure at spatial infinity,  $C_P=\gamma C_V$, $C_V = R/(\gamma-1)$, $\gamma>1$ and $R>0$ are two constants. 
Given relation \eqref{press}, the set of unknowns for system \eqref{Lim} can be reduced to $\{\vr,\vu,\pi\}$ or to $\{\vu,\pi,T\}$  equivalently.\\
System \eqref{Lim} is complemented by the Newton rheological law for the viscous tensor, namely
$$\vc{S} = -2 \mu D(\vu) -\lambda \Div \vu \,  \vc{I}
$$
and the Fourier law for the heat flux
\eq{\vc{Q}=-k\nabla T,\label{Q_form}}
where $D(\vu) = \frac{1}{2}(\nabla \vu + \nabla^t \vu)$, $\vc{I}$ is the identity matrix, 
$\mu$ and $\lambda$ are the viscosity coefficients and $k$ is the heat conductivity coefficient, all assumed to depend smoothly on the temperature, thus, by relation \eqref{press} on the density.\\
 Using \eqref{press} and \eqref{Q_form}, system \eqref{Lim} may be replaced by
 \begin{equation}\label{Limref}
\begin{array}{c}
\partial_t \vr + \Div (\vr \vu) = 0, \\
\pt (\vr \vu) + \Div (\vr \vu \otimes \vu)
         +\nabla \pi = -\delta \Div\vc{S} ,\\
\Div \vu =\displaystyle \frac{(\gamma -1) }{\gamma R} \kappa \Div\lr{{k}\Grad \lr{ \frac{1}{\vr}}}.
\end{array}
\end{equation} 
 Aassuming $\delta =1$ {and including dependence on the temperature in the coefficient $k$}, the following system is obtained
 \begin{equation}\label{Limref2}
\begin{array}{c}
\partial_t \vr + \Div (\vr \vu) = 0, \\
\pt (\vr \vu) + \Div (\vr \vu \otimes \vu)
         +\nabla \pi =  2 \Div (\mu(\vr) D(\vu)) + \nabla (\lambda(\vr) \Div \vu) ,\\
\Div \vu = - 2  \kappa \Delta \varphi(\vr),
\end{array}
\end{equation} 
with $\varphi$--an increasing function of $\vr$.  
Depending on the
nonlinearity $\varphi$, such systems are used to model various phenomena like motion of mixtures and avalanches, salt and  pollutant spreading or combustion. 
In the recent paper \cite{GoVa} a more complex system is derived to model a flow of mixture in the multi-dimensional setting and investigated in the one-dimensional domain.
    
 \medskip
   
 Several authors, such as {\sc H. Beir\~{a}o Da Veiga} or {\sc P. Secchi}  \cite{Ve, Se2} have considered the problem of existence of local strong solutions to system \eqref{Limref2}.  
 The interested reader is referred to the recent interesting work by {\sc R. Danchin} and {\sc X. Liao} \cite{DaLi12} where the existence of global solutions in homogeneous Besov spaces with critical regularity is proven assuming the initial density close to a constant and the initial velocity
 small enough.  
 Concerning the global in time existence,    {\sc A.~Kazhikov, S.~Smagulov} showed in \cite{KaSm77} that system \eqref{Limref2} with modified convective term and small $\kappa$ possesses a global in time generalized solution that is unique in two-dimensional domain.  
     {\sc P. Secchi} in \cite{Se1} proved the existence of a (unique) global solution for two-dimensional flows  when the diffusion coefficient $\kappa$ is small. He also considered the convergence (as $\kappa \to 0$) towards the corresponding solutions of the  nonhomogeneous Navier--Stokes system for two- and three-dimensional case. In \cite{DaLi12}, the authors  proved the global existence result  in critical spaces if the density is close to a constant and if  the initial velocity is small enough.  In \cite{PLL} {\sc P.--L.~Lions} showed, in two-dimensional case,  that for a 
positive conductivity  coefficient and $\varphi= {-1/\vr}$, a small perturbation of a constant density provides a global existence of weak solutions without restriction on the initial velocity.  He has yet left
a generalization of his result to the three-dimensional case as an open problem.
     
 \medskip
 
  The first global existence result of weak solutions without smallness assumption was obtained by {\sc D. Bresch}, {\sc E.H. Essoufi} and {\sc M. Sy} in \cite{BrEsSy07} when a certain algebraic relation between $\mu$ and $\kappa \varphi$ is assumed, namely
\eq{\varphi'(s) =  \mu'(s)/s \quad \text{and} \quad \kappa=1,\label{change1}}
for which the third equation in \eqref{Limref2} becomes
$$\Div\vu=-2\Div(\vr^{-1}\Grad\mu(\vr)).$$
   Later on, {\sc X. Cai, L. Liao, Y. Sun} proved the uniqueness of this solution in the two-dimensional case \cite{CLS12}. Recently, this algebraic relation was also used by {\sc X. Liao}  \cite{X.Liao} to show existence of weak solutions and. In the two-dimensional case, she also showed uniqueness of this solution in the critical non-homogeneous Besov spaces.
 
 \medskip
    
 In the first part of this paper, we show how to relax relation \eqref{change1}. More precisely, we prove existence of global weak solutions assuming only
 \begin{equation}\label{rel1}
 \varphi'(s) =  \mu'(s)/s \quad\text{ and }  \quad 0 < \kappa <  1
\end{equation}
which implies
\eq{\Div\vu=- 2\kappa\Div(\vr^{-1}\Grad\mu(\vr)).\label{non_phi}}
This result may be viewed as  a generalization of the particular case $\kappa=1$ studied  by {\sc D.~Bresch} $\&$ {\it al.} in \cite{BrEsSy07}  to the case when $\kappa$ is any constant from the interval $(0,1)$. It is based on the estimate of new mathematical entropy \eqref{hypo},  whose prototype for $\kappa=1$ was proposed in \cite{BrEsSy07}. 
Another interesting point of this paper is a construction of approximate solutions. The general relation \eqref{rel1} leads to presence of higher order terms in the momentum equation and thus the complexity of construction is significantly higher than the one from \cite{BrEsSy07, X.Liao}. It uses an original augmented regularized system \eqref{regularized} of parabolic type.
 Let us also announce that ideas developed in this paper will be adapted to handle the case compressible  Navier--Stokes equations in the continuation of this paper (see \cite{BrDeZa}).

\medskip

In the second part of the paper we use the existence result and the uniform estimates for $\kappa\in(0,1)$ to prove convergence of weak solution  towards the  corresponding solutions of the non-homogeneous incompressible Navier--Stokes system  (as $\kappa$ tends to $0$)  and solutions of the {\it Kazhikhov-Samgulov}-type system  (as $\kappa$ tends $1$).\\
The existence of global weak solutions for the non-homogenous incompressible
Navier-Stokes equations ($\kappa=0$) has been investigated: by {\sc S. N. Antontsev, A. V. Kazhikhov and V. N. Monakhov} \cite{AnKaMo90} for a constant viscosity and initial density bounded and far away from vacuum, by R. Danchin and P.B. Mucha for discontinuous initial density \cite{DaMu}, by  {\sc J. Simon} \cite{Sim90}  for constant viscosity and bounded density with possible vacuum and by {\sc P.--L.~Lions} \cite{PLLI} for non-degenerate density dependent viscosity coefficients with $L^\infty(\Omega)$ bound for the density and possible vacuum.  Here,  the limit passage $\kappa\to0$ is performed in the case when the initial density is far away from zero and bounded.

\bigskip

The third part of the paper is devoted to study of a more general form of nonlinearity than \eqref{non_phi}.  We show how to relax the algebraic relation \eqref{rel1}, more precisely, we
  consider the following system
   \begin{equation}\label{Limref1}
\begin{array}{c}
\partial_t \vr + \Div (\vr \vu) = 0, \\
\pt (\vr \vu) + \Div (\vr \vu \otimes \vu)
         +\nabla \pi =  2 \Div (\mu(\vr) D(\vu)) + \nabla (\lambda(\vr) \Div \vu) ,\\
\Div \vu = -  \Delta \tilde \varphi(\vr),
\end{array}
\end{equation} 
with $\tilde \varphi$--an increasing function of $\vr$.  Then, defining a new function $\tilde\mu(s)$
such that $\tilde \mu'(s) = s \tilde \varphi'(s)$, we propose three inequalities relating 
$\mu$ and $\tilde \mu$  (see \eqref{c_gen}) which yield the global existence of weak solutions.   
 These inequalities allow to generalize  the mathematical entropy  introduced in the first part to \eqref{integ1}. 
 As a corollary of these inequalities, we
 prove that a small perturbation of a constant initial density provides a global existence of weak solutions. This was observed already  by  {\sc P.--L.~Lions} in \cite{PLL} for two-dimensional domain, but thanks to our observation it is true also in the  three-dimensional case with no smallness assumption on the initial velocity. Hence, our result sheds a new light on two open questions given in the book by {\sc P.--L.~Lions} \cite{PLL}. Firstly, we prove the global in time existence of solutions with no smallness assumption and without the severe restriction on the form of nonlinearity \eqref{change1} that we replace by a relaxed one namely \eqref{rel1}. Secondly, we provide an existence of global weak solutions with restriction on the initial density but {\it not on the initial velocity} for $d = 2$ and  $d=3$. 

\bigskip

  In the two last sections, we  present applications of our result to the gaseous mixture  and ghost effect systems. More precisely, in the Section \ref{S:mix}, we show how to extend the result by {\sc E. Embid} or  T. {\sc Alazard} (see \cite{Embid87} and \cite{Al1}) about local existence of strong solutions to the framework of weak solutions for two-component mixture model.
   In Section \ref{S:ghost}, we discuss the system derived in  \cite{LeSuTr} by {\sc C.D. Levermore, W. Sun,  K. Trivisa} as a low Mach number limit for classical solutions of the compressible Navier-Stokes equations with dispersive corrections. 
\begin{equation*}
\begin{array}{c}
\vr T=1, \qquad 
\partial_t \vr + \Div (\vr \vu) = 0, \\
\pt (\vr \vu) + \Div (\vr \vu \otimes \vu)
         +\nabla P^*= - \Div\vc{\Sigma}-\Div\tilde{\vc{\Sigma}} ,\\
\frac{5}{2}\Div\vu =\Div\lr{k(T)\Grad T},
\end{array}
\end{equation*} 
where $k(T)$ is the heat conductivity coefficient $k>0$, while the two parts of stress tensor $\vc{\Sigma}$ and $\tilde{\vc{\Sigma}}$ are defined as follows
\eqh{
\vc{\Sigma}&=\mu(T)\lr{\Grad\vu+\Grad^t\vu-\frac{2}{3}\Div\vu \vc{I}},\\
\tilde{\vc{\Sigma}}&=\tau_1(\vr,T)\lr{\Grad^2T-\frac{1}{3}\lap T\vc{I}}+\tau_2(\vr,T)\lr{\Grad T\otimes\Grad T-\frac{1}{3}|\Grad T|^2\vc{I}},
}
where $\tau_1,\tau_2$ are transport coefficients with $\tau_1>0$. For a specific choice of physical viscosity,  transport coefficient and heat-conductivity coefficient, we show how to apply our theoretical result to get global existence of weak solutions to a ghost effect system. To our knowledge, this gives a first answer to a question about global existence of weak solutions to such kinds of systems.  

 \medskip
 
 As the title suggests, the present paper is the first part of a series. We refer the interested reader to part II  \cite{BrDeZa} for the extension of the above results to the case of two-velocity hydrodynamics in compressible Navier-Stokes equations with degenerate viscosities: the construction of approximate solutions will follow the same lines than in the present paper but will be
fully given for reader's convenience.

\medskip
\noindent {\bf Remark about the notation}: In the sequel $c$ denotes generic positive constant (possibly large) that may change from line to line.

\section{Main results}
{This section is devoted to presentation of our main result. First, we state the global in time existence result in the case when a special algebraic relation between $\vp$ and $\mu$ is assumed. In the second part we relax this algebraic equality to inequality and we formulate  the second  existence result  covering more general form of $\vp$ and $\mu$.}

 \subsection{The case of $\varphi$ and $\mu$ related by \eqref{rel1}}

 \medskip
 
 \noindent {\bf Reformulation of the system.}   Before formulating our main result we rewrite system \eqref{Limref2} in a different form. These forms are equivalent provided solutions to \eqref{Limref2} are sufficiently regular.
   Let us first introduce the following solenoidal vector field
\eq{&\vw=\vu+ 2\kappa \Grad \vp(\vr),\\
&\Div\vw=0.\label{def_w}}
Using this notation, system \eqref{Limref1} can be rewritten as
\begin{equation}\label{main1}
\begin{array}{c}
\pt\vr+ \Div(\vr \vu)  = 0,\\
\pt\lr{\vr \vw}+ \Div  (\vr \vu \otimes \vw) - 2 \Div (\mu(\vr) D(\vu)) + 2 \kappa \Div (\mu(\vr) \Grad^t \vu) + \Grad {\pi_1} =\vc{0},\\
\vw=\vu+2\kappa \Grad \vp(\vr),\\
\Div\vw = 0,
\end{array}
\end{equation}
where 
$$\pi_1= \pi + 2(\mu'(\vr)\vr-\mu(\vr))\Div\vu- \lambda(\vr)\Div\vu.$$ 
To see this it is enough to multiply the continuity equation by $\mu'(\vr)$ and write the corresponding equation for $\mu(\vr)$ 
\eq{\pt\mu(\vr) + \Div(\mu(\vr)\vu) + (\mu'(\vr)\vr-\mu(\vr))\Div\vu = 0.\label{mu}}
Differentiating it with respect to space and employing \eqref{rel1}, we obtain
\eq{\pt\lr{\vr\Grad\vp(\vr)} + \Div(\vr \vu\otimes \Grad \vp(\vr)) 
      + \Div(\mu(\vr) \Grad^t \vu) + \Grad\lr{(\mu'(\vr)\vr-\mu(\vr))\Div\vu}= \vc{0}.\label{vp}}
Thus the second equation of \eqref{main1} is obtained using definition of \eqref{def_w} and the momentum equation from \eqref{Limref1}. $\Box$

\medskip

\noindent {\bf Hypothesis, and definitions of weak solution.}
System \eqref{Limref2} is supplemented by the periodic boundary conditions, i.e. 
$$\Omega=\mathbb{T}^3.$$ 
and the initial conditions
\eq{\vr|_{t=0} =\vr^0,\quad \vu|_{t=0}=\vu^0, \quad \vw|_{t=0}=\vw^0=\vu^0+2\kappa\Grad\vp(\vr^0) \quad \text{in}\  \Omega.
\label{init}}
We assume that the initial conditions satisfy 
\eq{\label{init_cond}
\sqrt{(1-\kappa)\kappa}\vr^0 \in H^1(\Omega), \qquad 0< r\leq \vr^0\leq R<\infty, \qquad 
\vw^0 \in H,
}
uniformly with respect to $\kappa$, where 
\eqh{H=\{\vc{z}\in L^2(\Omega); \ \Div\vc{z}=0\}\quad \text{and} \quad V=\{\vc{z}\in W^{1,2}(\Omega);\  \Div\vc{z}=0\}.}
We look for generalized weak solutions of system \eqref{main1} in the sense of the following definition.
\begin{df}[Global weak solution in terms of $\vw$]\label{Def1}
The couple of functions $(\vr,\vw)$ is called a global weak solution to system \eqref{main1} and \eqref{init} 
if the following regularity properties are satisfied 
\begin{equation*}
\begin{gathered}
0<r\leq \vr\leq R<\infty, \quad a.e. \ in\ (0,T)\times\Omega,\\
\vr\in L^\infty(0,T;H^1(\Omega))\cap L^2(0,T;H^2(\Omega)),\\
\vw\in L^\infty(0,T;H)\cap L^2(0,T;V),
\end{gathered}
\end{equation*}
 the $\kappa$-entropy estimate holds
 \eq{\label{hypo}
 &\sup_{\tau\in [0,T]}  \intO{\vr\lr{\frac{|\vw|^2}{2}+ (1-\kappa)\kappa\frac{|2\Grad\vp(\vr)|^2}{2}}(\tau)}
   + 2\kappa\intTO{\mu(\vr) |A(\vw)|^2}\\
&+ 2(1-\kappa)\intTO{ \mu(\vr) \left|D(\vw-2\kappa\nabla\varphi(\vr))
    + \frac{2\kappa }{d}\lap \varphi(\vr)\,  \vc{I}\right|^2} \\
& +2(1-\kappa) \intTO{\lr{\frac{(1-d)}{d}\mu(\vr)+\mu'(\vr)\vr}|2\kappa\lap\vp(\vr)|^2}\\
&\hskip5cm  \le \intO{\vr\lr{\frac{|\vw|^2}{2}+ (1-\kappa)\kappa\frac{|2\Grad\vp(\vr)|^2}{2}}(0)}, 
}
where $A(\vw)=\frac{1}{2}(\Grad\vw-\Grad^t\vw)$ and the equations of system \eqref{main1} hold in the sense of distributions.\\
More precisely, the continuity equation
\eq{\intTO{\vr\pt\phi}+\intTO{\vr\vw\cdot\Grad\phi}- 2\kappa\intTO{\vr\Grad\vp(\vr)\cdot\Grad\phi}=-\intO{\vr^0\phi(0)}\label{weak_cont}}
holds for $\phi\in C^\infty([0,T]\times\Omega)$, s.t. $\phi(T)=0$.\\
 The momentum equation written in terms of $\vw$
\eq{
&\intTO{\vr\vw\cdot\pt\vphi_1}+\intO{\vr(\vw- 2\kappa\Grad\vp(\vr))\otimes\vw:\Grad\vphi_1}\\
&\quad-2\intTO{\mu(\vr)D(\vw):\Grad\vphi_1}+ 2\kappa\intTO{\mu(\vr)\Grad^t\vw:\Grad\vphi_1}\\
&\quad+ 4(1-\kappa)\kappa\intO{\mu(\vr)\Grad \Grad\vp(\vr):\Grad\vphi_1} =-\intO{\vr^0\vw^0\cdot\vphi_1(0)}\label{weak_mom}}
holds for $\vphi_1\in (C^\infty([0,T]\times\Omega))^3$, s.t.  $\Div\, \vphi_1=0$ and $\vphi_1(T)=\vc{0}$.\\
The equation for $\nabla\varphi(\vr)$ 
\eq{
&\intTO{\vr\Grad\vp(\vr)\cdot\pt\vphi_2}+\intO{\vr(\vw- 2\kappa\Grad\vp(\vr))\otimes\Grad\vp(\vr):\Grad\vphi_2}\\
&\quad - \kappa \intTO{\mu(\vr)\Grad^2\vp(\vr) : \Grad\vphi_2}
     - \kappa\intTO{(\mu'(\vr)\vr -\mu(\vr)) \Delta \vp(\vr) \Div \vphi_2}\\
&\quad +  \intO{ \mu(\vr) \nabla^{t} \vw :\Grad\vphi_2}
 =-\intO{\vr^0\Grad\vp(\vr^0)\cdot\vphi_2(0)}\label{weak_graphi}}
holds for $\vphi_2\in (C^\infty([0,T]\times\Omega))^3$, s.t. $\vphi_2(T)=0$.
\end{df}
Let us now specify the notion of a weak solution to the original system \eqref{Limref2}. 
\begin{df}[Global weak solution in terms of $\vu$]\label{Def2}
The couple $(\vr,\vu)$ is called a weak solution to system \eqref{Limref2} and \eqref{init}
if the following regularity properties are satisfied 
\begin{equation*}
\begin{gathered}
0<r\leq \vr\leq R<\infty, \quad a.e. \ in\ (0,T)\times\Omega,\\
\vr\in L^\infty(0,T;H^1(\Omega))\cap L^2(0,T;H^2(\Omega)),\\
\vu\in L^\infty(0,T;L^2(\Omega))\cap L^2(0,T;W^{1,2}(\Omega)),
\end{gathered}
\end{equation*}
and the equations of system \eqref{Limref2} hold in the sense of distributions.\\
More precisely, the mass equation is satisfied in the following sense 
\eq{\intTO{\vr\pt\phi}+\intTO{\vr\vu\cdot\Grad\phi}=-\intO{\vr^0\phi(0)}\label{weak_contu}}
 for $\phi\in C^\infty([0,T]\times\Omega)$, s.t. $\phi(T)=0$. \\
 The momentum equation is satisfied in the following sense
\eq{ \intTO{\vr\vu\cdot\pt\vphi}+\intO{\vr \vu \otimes\vu:\Grad\vphi}
    -2\intTO{\mu(\vr)D(\vu):\Grad\vphi}\\
     =-\intO{\vr^0\vu^0\cdot\vphi(0)}\label{weak_momu}}
for $\vphi\in (C^\infty([0,T]\times\Omega))^3$ s.t.  $\Div\vphi=0$ and $\vphi(T)=\vc{0}$.\\
 The constraint 
 $$\Div \vu = - 2\kappa \Delta\varphi(\vr)$$
  is satisfied in $L^2(0,T;L^2(\Omega))$.
\end{df}
Defining $\vu = \vw - 2\kappa \Grad\vp(\vr)$, we see that the weak solution from Definition \ref{Def1} gives a weak solution from Definition \ref{Def2}. Indeed, using the definition of $\vu$ and the fact that $\vw$ is a divergence free vector field, the weak formulation of the momentum equation from Definition \eqref{Def2} is obtained by choosing $\vphi_1 = \vphi_2 = \vphi$ in \eqref{weak_mom} and \eqref{weak_graphi}, multiplying the second equation by $\kappa$ and subtracting it from the first one.

\bigskip

\noindent{\bf Formulation of the main results.} The first main result of this paper concerns the global in time existence of weak solutions to system \eqref{main1}.

\begin{thm}\label{T_main} Let $0<\kappa < 1$, $\mu$ be increasing function of class $C^1([r,R])$ such that $\mu\geq \Un{c}>0$ on $[r,R]$, assume \eqref{rel1} and let $\mu(\vr)$ satisfy the following condition on $[r,R]$
   \eq{\label{important}
    \lr{\frac{1-d}{d}\mu(\vr)+\mu'(\vr)\vr}\geq \Un{c}_1> 0.
   }
If the initial data $(\vr^0,\vw^0)$ satisfy  \eqref{init_cond}, then there exists at least one global weak solution  $(\vr,\vw)$ of system \eqref{main1}, in the sense of Definition \ref{Def1}.\\
Moreover, this solution satisfies the following estimates
%
\eqh{\sqrt\kappa\|\vr\|_{L^2(0,T;H^1(\Omega))}+\|\vw\|_{L^\infty(0,T;H)}+\sqrt{\kappa}\|\vw\|_{L^2(0,T; V)}\leq c,}
\eqh{\sqrt{(1-\kappa)\kappa}\|\vr\|_{L^\infty(0,T;H^1(\Omega))}+
    \sqrt{(1-\kappa)}\kappa\|\vr \|_{L^2(0,T;H^2(\Omega))}
+\sqrt{(1-\kappa)} \|D(\vu)\|_{L^2(0,T;L^2(\Omega))}\leq c,}
uniformly with respect to $\kappa$.\\
If, in addition to \eqref{init_cond}, $\kappa \vr^0\in{H^1(\Omega)}$ uniformly with respect to $\kappa$,  then 
\eq{\label{h_reg_rho}
\kappa^{3/2}\|\pt\mu(\vr)\|_{L^2(0,T; L^{3\over2}(\Omega))}+
  \kappa\| \mu(\vr)\|_{L^\infty(0,T; H^1(\Omega))}+
  \kappa^{3/2}\| \mu(\vr)\|_{L^2(0,T; H^2(\Omega))}\leq c
  }
and 
$$\sqrt \kappa\| \Grad \vu\|_{ L^2(0,T;L^2(\Omega))}\leq c$$
uniformly with respect to $\kappa$.
\end{thm}
Regularity of weak solutions follows from the $\kappa$-entropy inequality \eqref{hypo}.
It is important to note that 
$$|\vw|^2 + (1-\kappa)\kappa|2\nabla\varphi(\vr)|^2 =
    (1-\kappa)|\vu|^2 + \kappa  |\vu+2\nabla\varphi(\vr)|^2,$$
 therefore \eqref{hypo} reflects a two-velocity hydrodynamics with 
 a mixing ratio $\kappa$ in the spirit of the works by {\sc S.M. Shugrin} (see for instance {\rm \cite{Sh}})
 but with an incompressible mean velocity $\vw$.
 The interested reader is referred to the second part of this series of papers {\rm \cite{BrDeZa}}  where  we make a link between the compressible two-velocity hydrodynamics and the compressible-incompressible two-velocity hydrodynamics studied here.

\medskip

Having proven this theorem we investigate the limit $\kappa\to0$ to recover the usual non-homogenous incompressible Navier-Stokes equations and the limit $\kappa\to 1$ to recover the so-called Kazhikhov-Smagulov system. We have the following theorem

\begin{thm}\label{T_main2}
   Let $(\vr_\kappa,\vw_\kappa)$ be the sequence of global generalized weak solutions  obtained in the previous theorem, emanating from the initial data $(\vr_\kappa^0,\vw_\kappa^0)$. Assume that the initial data satisfy \eqref{init_cond} and
\begin{equation*}
\begin{gathered}
\sqrt{(1-\kappa)\kappa} \vr^0_\kappa\quad\text{ uniformly bounded
in }  H^1(\Omega),\\
\vr^0_\kappa\to\vr^0\quad\text{strongly in } L^p(\Omega) \hbox{ for all } 1< p<+\infty,\\
     \vw^0_\kappa= \vu^0_\kappa + 2\kappa \nabla \varphi(\vr^0)  \to \vw^0 
        \quad \hbox{ strongly in } V.
\end{gathered}
\end{equation*}
Then we have two cases:
\begin{itemize}
\item when $\kappa\to 0$ then there exists a subsequence, denoted again by $\kappa$, and limit functions $(\vr,\vu)$ such that 
\begin{equation*}
\begin{gathered}
\vr_\kappa\to\vr \hbox{ in }  C([0,T]; L^p(\Omega)) \hbox{ for all } 1< p<+\infty,\\
\vw_\kappa \to \vu \quad \text{weakly in } L^2(0,T;V) \text{ and weakly$^*$ in } 
     L^\infty(0,T;H)
\end{gathered}
\end{equation*}
and $(\vr,\vu)$ satisfies the following system
\begin{equation}\label{mainINS}
\begin{gathered}
\pt\vr+ \Div(\vr \vu)  = 0,\\
\pt\lr{\vr \vu}+ \Div  (\vr \vu \otimes \vu) - 2 \Div (\mu(\vr) D(\vu)) + \Grad {\pi_1} =\vc{0},\\
\Div\vu = 0,
\end{gathered}
\end{equation}
 in the sense specified in \eqref{weak_contu} and \eqref{weak_momu}, where $D(\vu)=\frac{1}{2}(\Grad\vu+\Grad^t\vu)$.

\item when $\kappa\to 1$, assuming in addition to \eqref{init_cond} that $\kappa \vr^0\in{H^1(\Omega)}$,  then there exists a subsequence, demoted again by $\kappa$, and limit functions $(\vr,\vu)$ such that 
\begin{equation*}
\begin{gathered}
\vr_\kappa\to\vr \hbox{ in }  C([0,T];L^p(\Omega))  \hbox{ for all } 1< p<+\infty, \\
(\kappa \mu'(\vr_\kappa))^{1/2} \nabla \vr_\kappa \to
(\mu'(\vr))^{1/2} \nabla \vr \hbox{ in } L^2(0,T;L^2(\Omega)),\\
\vw_\kappa\to \vw  \quad \text{weakly in } L^2(0,T;V) \text{ and weakly$^*$ in } 
     L^\infty(0,T;H)
\end{gathered}
\end{equation*}
and $(\vr,\vu)$ satisfies the following system
\begin{equation}\label{main2}
\begin{gathered}
\pt\vr+ \Div(\vr \vu) = 0,\\
\pt\lr{\vr \vw}+ \Div  (\vr \vu \otimes \vw) - 2 \Div (\mu(\vr) A(\vw)) + \Grad {\pi_1} =\vc{0},\\
\vw = \vu + 2 \nabla\varphi(\vr),  \qquad
\Div\vw = 0,
\end{gathered}
\end{equation}
in the sense specified in Definition \ref{Def1}, where $A(\vu)=\frac{1}{2}(\Grad\vu-\Grad^t\vu)$.

\end{itemize}
\end{thm}

\begin{rmk}
Observe that the following compatibility condition is satisfied for $\kappa=0$
$$\vw^0=\vu^0 \hbox{ with } \Div \vu^0 = 0.$$
\end{rmk}

\begin{rmk}
Note that for $\kappa=1$, we recover the existence result by {\sc D. Bresch} {\it et al.} in  {\rm \cite{BrEsSy07}}.
\end{rmk}

 \subsection{The case of general $\vp$ and $\mu$}
  In this section, we relax the algebraic relation between $\varphi$ and $\mu$
   \eqref{rel1}. More precisely, the diffusion equation 
$$\Div \vu = - 2 \kappa \Delta\varphi(\vr) \qquad \text{ with} \qquad
    \varphi'(s)= \mu'(s)/s$$ 
is replaced by
\eq{\Div \vu = -2 \Delta \widetilde\varphi(\vr)\qquad \text{ with} \qquad
    \widetilde\varphi'(s)= \widetilde\mu'(s)/s\label{gen_pm}}
with some new function $\widetilde \mu(s)$ which  is related to $\mu(s)$  only by some inequality. \\

  We prove the following result
 
 \begin{thm}\label{T_main3}  Let $\mu$ be an increasing function of class $C^1([r,R])$ such
 that $\mu\ge c>0$ on $[r,R]$ and $\tilde \varphi$ an increasing function of class $C^1([r,R])$.
 Assume that the initial data satisfy $\vr^0\in H^1(\Omega)$, $\vu^0\in (L^2(\Omega))^3$  with
 $$ 0< r\le \vr^0\le R < +\infty, \qquad 
      {\Div} (\vu^0 + 2\nabla \tilde \varphi(\vr^0)) = 0.$$
Moreover assume that there exist positive constants $c,c',\xi$ such that

\begin{equation}
\begin{gathered}
c\leq\min_{\vr\in[r,R]}\lr{\mu(\vr) -\tilde \mu(\vr)},\\
c'\leq\tilde\mu'(\vr)\vr+\frac{1-d}{d}\tilde\mu(\vr),\\
\max_{\vr\in[r,R]} \frac{(\mu(\vr)-\tilde\mu(\vr)-\xi\tilde\mu(\vr))^2}{ 2\lr{\mu(\vr) -\tilde \mu(\vr)}}
\leq \xi 
\min_{\vr\in[r,R]}\lr{\tilde\mu'(\vr)\vr+\frac{1-d}{d}\tilde\mu(\vr)}.
\end{gathered}
 \label{c_gen}
 \end{equation}
Then there exists a global weak solution to System \eqref{Limref1}.
\end{thm}

For this theorem, we will only prove the estimates because the construction and stability
process follows the lines given for Theorem \ref{T_main}. For the sake of completeness,
we give an example which show that this theorem gives an answer to a question
formulated by {\sc P.--L. Lions} in \cite{PLL}, Chapter 8.8 (see Proposition \ref{Example}). 
Global existence of weak solution for initial density close to constant and large velocity field in  the
three-dimensional in space case.  

 \medskip
 
 \section{Proof of Theorem \ref{T_main}}
This section is dedicated to the proof of existence of weak solutions to system \eqref{main1}. We first assume that $\vr,\vw$ are smooth enough and present the derivation of a-priori estimates necessary to understand the main idea of construction of solution. The latter is presented in Section \ref{SS:Constr}.
 
 \subsection{A priori estimates}
\noindent {\bf Maxiumum principle and $H^1$ bounds on the density.}
   First, applying the standard maximum principle for the continuity equation
$$\pt \vr + \vw\cdot\Grad \vr - 2 \kappa \lap\mu(\vr) = 0$$
 we deduce that 
\eq{0<r\leq\vr \leq R<\infty\label{max_vr},}
and the basic energy estimate gives
\eqh{\|\vr\|_{L^\infty(0,T;L^2(\Omega))}+ \sqrt\kappa \|\vr\|_{L^2(0,T; H^1(\Omega)}\leq c}
uniformly with respect to $\kappa$.

 \bigskip
 
 \noindent {\it The $\kappa$-entropy and its consequences.}
Our next goal is to derive an original estimate on $(\vr,\vw,\Grad\varphi(\vr))$. More precisely we 
prove that the following mathematical entropy holds
\begin{prop} Let $(\vr,\vw)$ be sufficiently smooth solution to \eqref{main1}, then $(\vr,\vw)$ satisfy inequality\eqref{hypo}.
\end{prop} 

\noindent {\it Proof.} We rewrite equation for $\vw$ from \eqref{main1} in slightly different form
\eq{\label{eqw}
\pt\lr{\vr \vw}+ \Div  (\vr \vu \otimes \vw) - 2(1-\kappa) \Div(\mu(\vr) D(\vu)) - 2\kappa \Div (\mu(\vr) A(\vu)) + \Grad {\pi_1} =\vc{0},}
where $A(\vu)=\frac{1}{2}(\Grad\vu-\Grad^t\vu)$.
Multiplying this equation by $\vw$ and integrating by parts with respect to $\Omega$ we obtain
\eq{\label{a}
\frac{1}{2}\Dt\intO{\vr|\vw|^2}+ 2(1-\kappa)\intO{\mu(\vr) |D(\vu)|^2}+ 2\kappa\intO{\mu(\vr) |A(\vu)|^2}\\
+ 4(1-\kappa)\kappa\intO{\mu(\vr)\Grad\vu:\Grad^2\vp(\vr)}=0.
}
   In the case $\kappa= 1$, the high-derivative of $\varphi(\vr)$ does not appear in this equality
and global existence of weak solutions may be obtained in a simpler way. Such calculations were firstly performed in \cite{BrDeZa} and may be also found in \cite{X.Liao}.
   Here, however, we do not assume that $\kappa=1$, thus it is not yet sure that \eqref{a} provides uniform estimates. The problem is the presence of the last term $\mu(\vr)\Grad\vu:\Grad^2\vp(\vr)$ which does not have a sign.  To eliminate it we take advantage of  the renormalized continuity equation \eqref{mu}. We first test  its gradient \eqref{vp} by $\Grad\vp(\vr)$ and we integrate it with respect to $\Omega$ to get
\eq{\label{b}
\Dt\intO{\vr\frac{|\Grad\vp(\vr)|^2}{2}}-\intO{\mu(\vr)\Grad\vu:\Grad^2\vp(\vr)}-\intO{(\mu'(\vr)\vr-\mu(\vr))\Div\vu\lap\vp(\vr)}=0.
}
Next, we multiply \eqref{b} by $4(1-\kappa)\kappa$ and add it to \eqref{a}, we obtain
\eq{\label{aa}
\Dt\intO{\vr\lr{\frac{|\vw|^2}{2}+(1-\kappa)\kappa\frac{|2\Grad\vp(\vr)|^2}{2}}}
&+ 2(1-\kappa)\intO{\mu(\vr) |D(\vu)|^2}+ 2\kappa\intO{\mu(\vr) |A(\vu)|^2}\\
&+ 2(1-\kappa)\intO{(\mu'(\vr)\vr-\mu(\vr))|2\kappa \lap\vp(\vr)|^2}=0.
} 

To find an optimal condition for $\mu(\vr)$ we apply the following equality
\eq{
     \intO{ \mu(\vr)|D(\vu)|^2 }
     = \intO{ \mu(\vr) \left|D(\vu) - \frac{1}{d}\Div  \vu\,  \vc{I}\right|^2} 
     + \intO{  \frac{\mu(\vr)}{d} |\Div  \vu|^2}
   \label{Ddiv}}
 to the second integral in \eqref{aa}, we get \eqref{hypo}.
Thus,  for the l.h.s. to be nonnegative and to control the second gradient of $\vr$ independently from $\kappa$ one needs to assume
   \eqh{
    \lr{\frac{1-d}{d}\mu(\vr)+\mu'(\vr)\vr}\geq c > 0,
   }
which justifies the assumption \eqref{important} from Theorem \ref{T_main}. $\Box$\\
Equality \eqref{hypo}  will be generalized to the case of compressible Navier-Stokes system in \cite{BrDeZa}.

\begin{rmk}
For example, for $\mu(\vr)=\vr^\alpha$, we have $\alpha >  1-\frac{1}{d}$.
\end{rmk}
Integrating \eqref{hypo} with respect to time we get
\eqh{
&\intO{\vr\lr{\frac{|\vw|^2}{2}+ (1-\kappa)\kappa\frac{|2\Grad\vp(\vr)|^2}{2}}(T)}
    + 2\kappa\intTO{\mu(\vr) |A(\vu)|^2}\\
&+ 2(1-\kappa)\left[\intTO{ \mu(\vr) \left|D(\vu) - \frac{1}{d}\Div  \vu\,  \vc{I}\right|^2}
+ \intTO{\lr{\frac{(1-d)}{d}\mu(\vr)+\mu'(\vr)\vr}|\Div\vu|^2}\right]
\\ & \le \intO{\vr\lr{\frac{|\vw|^2}{2}+  (1-\kappa)\kappa\frac{|2\Grad\vp(\vr)|^2}{2}}(0)}.
}
Using \eqref{max_vr}, the above inequality can be used to deduce the following estimates
\eq{\|\vw\|_{L^\infty(0,T;H)}+
&   \sqrt{(1-\kappa)\kappa} \left\|\vr\right\|_{L^\infty(0,T;H^1(\Omega))}
  +\sqrt{\kappa}\|\vw\|_{L^2(0,T; V)} \\
& \hskip2cm + \sqrt{(1-\kappa)}\|D(\vu)\|_{L^2(0,T; L^2(\Omega))}
  + \sqrt{(1-\kappa)} \kappa \|\vr \|_{L^2(0,T; H^2(\Omega))}  \leq c.\label{v2}}
  
\bigskip

\noindent {\bf  Additional estimate on $\vr$.}   This estimate will be used in order to let $\kappa\to1$ in the proof of second part of Theorem \ref{T_main2}. Assuming that $\kappa \vr^0 \in H^1(\Omega)$ we show that it is possible to control 
$\kappa^{3/2}\lap\mu(\vr)$.  To see this let us first rewrite \eqref{mu}  as
\eq{\pt\mu(\vr) + \vw\cdot\Grad\mu(\vr)- 2\kappa\mu'(\vr)\lap\mu(\vr) = 0.\label{mu_w}}
Multiplying the above equation by $-\kappa^2\lap\mu(\vr)$ and integrating by parts we obtain
\eqh{\Dt\intO{\kappa^2|\Grad\mu(\vr)|^2}+2\intO{\kappa^3\mu'(\vr)|\lap\mu(\vr)|^2}\leq \intO{\kappa^2|\Grad\vw||\Grad\mu(\vr)|^2}.}
Note that for any $f\in H^2(\Omega)\cap L^\infty(\Omega)$ we have
\eq{\|\Grad f\|^2_{L^4(\Omega)}\leq c\|\lap f\|_{L^2(\Omega)}\|f\|_{L^\infty(\Omega)}.\label{G-N}}
Therefore, applying this inequality with $f=\mu(\vr)$ and using the uniform bound for $\sqrt\kappa \Grad \vw$ following from \eqref{v2},  we obtain
\eq{ \label{est}
\kappa \|\Grad\mu(\vr)\|_{L^\infty(0,T; L^2(\Omega))}+\kappa^{3/2}\|\lap\mu(\vr)\|_{L^2((0,T)\times\Omega)}\leq c\lr{\kappa\|\vr^0\|_{H^1(\Omega)},\sqrt{\kappa}\|\vw\|_{L^2(0,T;V)},R},
}
and thus, coming back to \eqref{mu_w} we easily check \eqref{h_reg_rho}.

\subsection{Construction of solution}\label{SS:Constr}

Construction of sufficiently smooth approximation is different when $0<\kappa <1$ than for the case $\kappa=1$. The latter  can found for instance in \cite{X.Liao} (adaptation of the usual process for example from \cite{PLL}). The difference between these two cases is due to the third term in \eqref{eqw} which does not vanish and which generates higher order term with respect to $\vr$. Here we have to consider the following augmented system with three unknowns $(\vr,\vw,\vv)$:
\begin{equation}\label{regularized}
\begin{array}{c}
\vspace{0.2cm}
\pt\vr+ \Div(\vr \vw) - 2 \kappa\lap\mu( \vr) = 0,\\
\hskip-8cm \pt\lr{\vr \vw}+ \Div  ((\vr \vw - 2 \kappa  \nabla\mu( \vr)) \otimes \vw)   \\
\vspace{0.2cm}
 \hskip2cm       - 2(1-\kappa)  \Div (\mu(\vr) D(\vw)) - 2\kappa \Div (\mu(\vr) A(\vw))  
   + \Grad {\pi_1} = - 2 \kappa (1-\kappa) \Div(\mu(\vr) \nabla\vv),\\
\hskip-8cm \partial_t(\vr \vv) + \Div((\vr \vw- 2 \kappa \nabla \mu( \vr))\otimes \vv) \\
\vspace{0.2cm}
\hskip2cm -2  \kappa \Div(\mu(\vr)\nabla\vv) - 2 \kappa \nabla ((\mu'(\vr)\vr - \mu(\vr))\Div \vv)
= - 2\Div(\mu(\vr) \nabla^{t} \vw), \\
\Div\vw = 0.
\end{array}
\end{equation}
Of course, we will have to prove that $\vv = 2 \nabla\varphi(\vr)$ to solve the initial system.
\bigskip
\bigskip

\noindent{\bf Full approximation.} Below we present the basic level of approximation procedure.

\bigskip

1. The continuity equation is replaced by its regularized version
\eq{&\pt\vr+ \Div(\vr [\vw]_\delta) - 2\kappa\, \Div([\mu'({\vr})]_\alpha \nabla \vr) = 0,\\
&\qquad\qquad \vr(0,x)=[\vr^0]_{\delta},
\label{A_cont}}
where $\alpha,\delta$ denote the standard regularizations with respect to time and space.
Such double index regularization may be found for instance in \cite{X.Liao}.

\bigskip

2. The momentum equation is replaced by its Faedo-Galerkin approximation with additional regularizing term $\ep\lr{\lap^2\vw-\Div((1+|\Grad\vw|^2)\Grad\vw)}$
\begin{equation}\label{FG}
\begin{split}
&\intO{\vr\vw(\tau)\cdot\vcg{\phi}}
-\inttauO{((\vr [\vw]_\delta -2 \kappa[\mu'({\vr})]_\alpha \nabla \vr) \otimes \vw):\Grad \vcg{\phi}}\\
&\quad + 2(1-\kappa)\inttauO{\mu(\vr)D(\vw):\Grad \vcg{\phi}}
+\ep\inttauO{\lr{\lap\vw\cdot\lap\vcg{\phi}+(1+|\Grad\vw|^2) \Grad\vw:\Grad\vcg{\phi}}}\\
&\quad+2 \kappa\inttauO{\mu(\vr)A(\vw): \Grad\vcg{\phi}}- 2\kappa(1-\kappa)\inttauO{\mu(\vr)\Grad\vv:\Grad\vcg{\phi}}\\
&=\intO{(\vr\vw)^{0}\cdot\vcg{\phi}},
\end{split}
\end{equation}
satisfied for any $\tau\in[0,T]$ and any test function $\vcg{\phi}\in X_{n}$, where $X_{n}=\operatorname{span}\{\vcg{\phi}_{i}\}_{i=1}^{n}$  and $\{\vcg{\phi}_{i}\}_{i=1}^{\infty}$ is an orthonormal basis in $V$ such that 
$\vcg{\phi}_{i} \in (C^\infty(\Omega))^3$ with $\Div\,  \vcg{\phi}_{i}=0$ for all $i\in \mathbb{N}$.

\medskip

3. The Faedo-Galerkin approximation for the artificial equation
\eq{\label{FGtheta}
&\intO{\vr\vv(\tau){\bf{\vcg{\xi}}}}
-\inttauO{((\vr [\vw]_\delta - 2 \kappa [\mu'({\vr})]_\alpha \nabla \vr) \otimes \vv):\Grad \vcg{\xi}}
+ 2\kappa\inttauO{\mu(\vr)\Grad\vv:\Grad \vcg{\xi}}\\
&\quad+ 2\kappa\inttauO{(\mu'(\vr)\vr - \mu(\vr))\Div \vv\ \Div\vcg{\xi}}- 2\inttauO{\mu(\vr)\Grad^t\vw:\Grad\vcg{\xi}}=\intO{(\vr\vv)^{0}\vcg{\xi}},
}
satisfied for any $\tau\in[0,T]$ and any test function $\vcg{\xi}\in Y_{n}$, where $Y_{n}=\operatorname{span}\{\vcg{\xi}_{i}\}_{i=1}^{n}$   and $\{\vcg{\xi}_{i}\}_{i=1}^{\infty}$ is an orthonormal basis in $W^{1,2}(\Omega)$
such that $\vcg{\xi}_{i}\in (C^\infty(\Omega))^3$ for all $i\in \mathbb{N}$.

\bigskip

\noindent{\bf Existence of solutions to the continuity equation.} For fixed $\vw\in C([0,T], X_n)$  we solve the continuity equation, which is now quasi-linear parabolic equation with smooth coefficients.
Thus, application of classical existence theory of {\sc Lady{\v{z}}enskaja, Solonnikov and Uralceva} \cite{LSU} (see for example Theorem 10.24 from \cite{FN}, which is a combination of Theorems 7.2, 7.3 and 7.4 from \cite{LSU}) yields the following result
\begin{thm}\label{LSU} Let $\nu \in (0,1)$. Suppose that the initial condition $\vr_\delta^0\in C^{2+\nu}(\Ov\Omega)$ is such that $0<r\leq \vr_\delta^0\leq R$ and satisfies the periodic boundary conditions. Then  problem \eqref{A_cont} possesses a unique classical solution $\vr$ from the class 
\begin{equation}\label{regvr}
V_{[0,T]}=\left\{\begin{array}{rl}
\vr&\in C([0,T];C^{2+\nu}(\Omega))\cap C^1([0,T]\times\Omega),\\
\pt\vr&\in C^{\nu/2}\left([0,T];C(\Omega)\right)
\end{array}\right\}
\end{equation}
and satisfying classical maximum principle 
\eq{0< r\leq\vr(t,x)\leq R. \label{max_gal}}
Moreover, the mapping $\vw\mapsto\vr(\vw)$ maps bounded sets in $C([0,T];X_{n})$ into bounded sets in $V_{[0,T]}$ and is continuous  with values in $C\big([0,T];C^{2+\nu'}(\Omega)\big)$, $0<\nu'<\nu<1$.

\end{thm}

\bigskip

\noindent{\bf Local existence of solutions to the Galerkin approximations.} Having a smooth solution to the continuity equation we may now solve the integral equations \eqref{FG} and \eqref{FGtheta} on possibly short time interval via fixed point argument.
More precisely, we want to prove that there exists  time $T=T(n)$ and $(\vw,\vv)\in C([0,T];X_{n})\times C([0,T]; Y_n)$ satisfying (\ref{FG}, \ref{FGtheta}).
To this purpose let us rewrite this equations as a fixed point problem
\eq{\label{T}
(\vw(t),\vv(t))&=\lr{ {\mathcal{M}}_{\vr(t)}\left[P_{X_n}(\vr\vw)^0+\int_{0}^{t}{{\mathcal{K}}(\vw)(s) {\rm d}s}\right],\ \mathcal{N}_{\vr(t)}\left[P_{Y_n}(\vr\vv)^0+\int_{0}^{t}{{\mathcal{L}}(\vv)(s) {\rm d}s}\right]}\\
&={\mathcal{T}}[\vw,\vv](t)
}
where $\vr=\vr(\vw)$ is a solution to the continuity equation as explained above,
$${\cal{M}}_{\vr(t)}:X_{n}\rightarrow X_{n},\quad\intO{\vr{\cal{M}}_{\vr(t)}[\vcg{\phi}]\cdot\vcg{\psi}}=\langle\vcg{\phi},\vcg{\psi}\rangle, \quad\vcg{\phi},\vcg{\psi}\in X_{n},$$
$${\cal{N}}_{\vr(t)}:Y_{n}\rightarrow Y_{n},\quad\intO{\vr{\cal{N}}_{\vr(t)}[\vcg{\xi}]\cdot\vcg{\zeta}}=\langle\vcg{\xi},\vcg{\zeta}\rangle, \quad\vcg{\xi},\vcg{\zeta}\in Y_{n},$$
$P_{X_n}$, $P_{Y_n}$ denote the projections of $L^2(\Omega)$ on $X_n$, $Y_n$, respectively, and ${\cal K}(\vw), {\cal L}(\vv)$ are the operators defined as
\eqh{
{\cal K}: &\quad X_n\to X_n,\\
\langle{\cal K}(\vw),\vcg{\phi}\rangle=&\intO{((\vr [\vw]_\delta - 2 \kappa [\mu'({\vr})]_\alpha \nabla \vr) \otimes \vw):\Grad \vcg{\phi}}
- 2(1-\kappa)\intO{\mu(\vr)D(\vw):\Grad \vcg{\phi}}\\
&-2 \kappa\intO{\mu(\vr)A(\vw): \Grad\vcg{\phi}}
   + 2\kappa(1-\kappa)\intO{\mu(\vr)\Grad\vv:\Grad\vcg{\phi}}\\
&-\ep\intOB{\lap\vw\cdot\lap\vcg{\phi}+(1+|\Grad\vw|^2) \Grad\vw:\Grad\vcg{\phi}},
}
\eqh{
{\cal L}: &\quad Y_n\to Y_n,\\
\langle{\cal L}(\vv),\vcg{\xi}\rangle=
&\intO{((\vr [\vw]_\delta - 2 \kappa [\mu'({\vr})]_\alpha \nabla \vr) \otimes \vv):\Grad \vcg{\xi}}
- 2\kappa\intO{\mu(\vr)\Grad\vv:\Grad \vcg{\xi}}\\
&-2 \kappa\intO{(\mu'(\vr)\vr - \mu(\vr))\Div \vv\ \Div\vcg{\xi}}+2\intO{\mu(\vr)\Grad^t\vw:\Grad\vcg{\xi}}.}
First let us observe that since $\vr(t,x)$ is strictly positive we have
\eqh{\|{\cal M}_{\vr(t)}\|_{L(X_n,X_n)},\ \|{\cal N}_{\vr(t)}\|_{L(Y_n,Y_n)}\leq\frac{1}{r}.}
Moreover
\eq{\|{\cal M}_{\vr^1(t)}-{\cal M}_{\vr^2(t)}\|_{L(X_n,X_n)}
+ \|{\cal N}_{\vr^1(t)}-{\cal N}_{\vr^1(t)}\|_{L(Y_n,Y_n)}\leq c(n,r^1,r^2)\|\vr^1-\vr^2\|_{L^1(\Omega)},\\
\label{XY}
}
 and by the equivalence of norms on the finite dimensional space we prove that
 \eq{\|{\cal K}(\vw)\|_{X_n}+ \|{\cal L}(\vv)\|_{Y_n}\leq c(n,r,R,\|\Grad\vr\|_{L^2(\Omega)},\|\vw\|_{X_n}, \|\vv\|_{Y_n} ).\label{KL}}
Next, we consider a ball ${\cal B}$ in the space $C([0,\tau];X_{n})\times C([0,\tau];Y_{n})$:
	$${\cal B}_{M,\tau}=\left\{(\vw,\vv)\in C([0,\tau];X_{n})\times C([0,\tau];Y_{n}):\|\vw\|_{C([0,\tau];X_n)}+\|\vv\|_{C([0,\tau];Y_n)}\leq M\right\}.$$
Using estimates \eqref{XY}, \eqref{KL}, \eqref{regvr} and \eqref{max_gal}, one can check that ${\cal T}$ is a continuous mapping of the ball ${\cal B}_{M,\tau}$ into itself and for sufficiently small $\tau=T(n)$ it is a contraction. Therefore, it possesses a unique fixed point which is a solution to \eqref{FG} and \eqref{FGtheta} for $T=T(n)$.

\bigskip

\noindent{\bf Global existence of solutions.} In order to extend the local in-time solution obtained above to the global in time one, we need to find uniform (in time) estimates, so that the above procedure can be iterated.
First let us note, that $\vw,\vv$ obtained in the previous paragraph have better regularity with respect to time. It follows by taking the time derivative of \eqref{T} and using the estimates \eqref{regvr}, \eqref{max_gal}, that
$$(\vw,\vv)\in C^1([0,\tau];X_{n})\times C^1([0,\tau];Y_{n}).$$
This is an important feature since now we can take time derivatives of \eqref{FG} and \eqref{FGtheta} and use the test functions $\vcg{\phi}=\vw$ and $\vcg{\xi}=\vv$, respectively. We then obtain
\begin{equation}\label{w1}
\begin{split}
&\Dt\intO{\vr\frac{|\vw|^2}{2}}
+ 2(1-\kappa)\intO{\mu(\vr)|D(\vw)|^2}+\ep\intOB{|\lap\vw|^2+(1+|\Grad\vw|^2) |\Grad\vw|^2}\\
&\quad+ 2\kappa\intO{\mu(\vr)|A(\vw)|^2}- 2 \kappa(1-\kappa)\intO{\mu(\vr)\Grad\vv:\Grad\vw}=0,
\end{split}
\end{equation}
and
\eq{\label{t1}
&\Dt\intO{\vr\frac{|\vv|^2}{2}}
+ 2\kappa\intTO{\mu(\vr)|\Grad\vv|^2}\\
&\quad+2 \kappa\intTO{(\mu'(\vr)\vr - \mu(\vr))(\Div \vv)^2}- 2\intTO{\mu(\vr)\Grad^t\vw:\Grad\vv}=0.
}
Therefore, multiplying the second equality by $(1-\kappa)\kappa$ and adding it to the first one, we obtain
\eq{\label{aa1}
\Dt&\intO{\vr\lr{\frac{|\vw|^2}{2}
+ (1-\kappa)\kappa\frac{|\vv|^2}{2}}}+ 2(1-\kappa)\intO{\mu(\vr) |D(\vw) -\kappa\nabla \vv)|^2} \\
&+\ep\intOB{|\lap\vw|^2+(1+|\Grad\vw|^2) |\Grad\vw|^2}\\
&+2 \kappa\intO{\mu(\vr) |A(\vw)|^2}+(1-\kappa)\intO{(\mu'(\vr)\vr-\mu(\vr))(\kappa\Div\vv)^2}=0.
}
Integrating the above estimate with respect to time, we obtain uniform estimate for $\vw$ and $\vv$ necessary 
to repeat the procedure described in the previous paragraph. Thus,  we obtain a global in time unique solution $(\vr,\vw,\vv)$ satisfying equations (\ref{A_cont}, \ref{FG}, \ref{FGtheta}).

\bigskip

\noindent{\bf Uniform estimates.} Below we  present uniform estimates that will allow us to pass to the limit with $\alpha$ and $n$ respectively.

First observe that multiplying continuity equation \eqref{A_cont} by $\vr_\alpha$ and integrating by parts with respect to $x$ gives
$$\frac{1}{2} \frac{d}{dt} \intO {\vr_\alpha^2} +2 \kappa \intO {[\mu'(\vr_\alpha)]_\alpha |\nabla\vr_\alpha|^2}  = 0.
$$
Integrating this equality with respect to time provides the following estimates 
\eq{\|\vr_\alpha\|_{L^\infty(0,T; L^2(\Omega))}+\|\sqrt{[\mu'(\vr_\alpha)]_\alpha }\nabla\vr_\alpha\|_{L^2(0,T;L^2(\Omega))}\leq c.\label{rho1}}
Moreover, the standard maximum principle gives boundedness of $\vr_\alpha$ from above and below. Indeed, multiplying equation \eqref{A_cont} by
$$\vr_\alpha^-=\max(0,r-\vr_\alpha)\quad\text{and}\quad \vr_\alpha^+=\min(0,R-\vr_\alpha),$$
respectively we obtain 
\eq{0< r\leq\vr_\alpha(t,x)\leq R. \label{max}}
\begin{rmk}
To prove these bounds one need to know that $\vr_\alpha\in L^2(0,T;W^{1,2}(\Omega))$, which doesn't follow from \eqref{rho1}. This problem could be solved by adding a small viscosity parameter $\alpha$ and considering $\widetilde{[\mu'(\vr)]_\alpha}=[\mu'(\vr)]_\alpha+\alpha$ in place of $[\mu'(\vr)]_\alpha$.
\end{rmk}
Next, using \eqref{max} and integrating \eqref{aa1} with respect to time we see that for $0<\kappa<1$ we have 
\eq{\|\vw\|_{L^\infty(0,T;H)}
+\|\vw_\alpha\|_{L^2(0,T;V\cap W^{2,2}(\Omega))}+
\|\Grad\vw_\alpha\|_{L^4(0,T;L^4(\Omega))}\\
+\|\vv_\alpha\|_{L^\infty(0,T;L^2(\Omega))}
+\|\vv_\alpha\|_{L^2(0,T;W^{1,2}(\Omega))}\leq c.
\label{ueu}}
\bigskip

\subsection{Passage to the limit with respect to $\alpha,n,\delta$ and $\varepsilon$}

\medskip

\noindent{\bf {Passage to the limit $\alpha\to0$.}} On the finite dimensional subspace all the norms are equivalent, therefore the space compactness of $\vw_\alpha$ and $\vv_\alpha$ is automatic. In fact, for $n$ fixed we also know that $\pt\vw_\alpha\in L^2(0,T; X_n)$, thus $\vw_\alpha\to\vw$ strongly in $L^2(0,T;X_n)$. The same can be deduced for $\vv_\alpha$.
The biggest problem is thus to pass to the limit in the term
\eq{[\mu'(\vr_\alpha)]_\alpha\Grad\vr_\alpha\otimes\vw_\alpha\label{conv_w0}}
which requires the strong convergence of the density and the weak convergence of the gradient of density.  Using \eqref{max} we can deduce that there exists $c>0$ and $\alpha_0$ such that for $\alpha<\alpha_0$ one has $[\mu'(\vr_\alpha)]_\alpha>c$ uniformly with respect to $\alpha$. Therefore  \eqref{rho1} implies that that up to a subsequence
\eqh{\vr_\alpha\to\vr\quad\text{weakly in }L^2(0,T; W^{1,2}(\Omega)),}
moreover $\pt\vr_\alpha\in L^2(0,T; W^{-1,2}(\Omega))$ and $\vr_\alpha \in L^\infty((0,T)\times\Omega)$, thus the Aubin-Lions lemma implies strong convergence of $\vr_\alpha$
\eqh{\vr_\alpha\to\vr\quad\text{strongly in }L^p((0,T)\times\Omega), \quad p<\infty.}
This justifies the passage to the limit in \eqref{conv_w0}. Therefore, one is able to pass to the limit $\alpha\to 0$ in both velocity equations \eqref{FG} and \eqref{FGtheta}. The limit functions $(\vr,\vw,\vv)=(\vr_n,\vw_n,\vv_n)$ satisfy the following system of equations:
\begin{itemize}
\item the momentum equation
\begin{equation}\label{FG_a}
\begin{split}
&\langle\pt\lr{\vr_n\vw_n}(t),\vcg{\phi}\rangle_{(X_n^*,X_n)}
-\intO{((\vr_n [\vw_n]_\delta - 2 \kappa \nabla \mu(\vr_n)) \otimes \vw_n)(t):\Grad \vcg{\phi}}\\
&\quad+ 2(1-\kappa)\intO{\mu(\vr_n)D(\vw_n)(t):\Grad \vcg{\phi}}
+ 2 \kappa\intO{\mu(\vr_n)A(\vw_n)(t): \Grad\vcg{\phi}}\\
&\quad-2\kappa(1-\kappa)\intO{\mu(\vr_n)\Grad\vv_n(t):\Grad\vcg{\phi}}
+\ep\intO{\lap\vw_n(t)\cdot\lap\vcg{\phi}}\\
&\quad+\ep\intO{(1+|\Grad\vw_n|^2) \Grad\vw_n(t):\Grad\vcg{\phi}}
=0,
\end{split}
\end{equation}
satisfied for $\vcg{\phi}\in X_n$, $t\in[0,T]$
\item the auxiliary equation for $\vv_n$
\eq{\label{FGtheta_a}
&\langle\pt\lr{\vr_n\vv_n}(t),\vcg{\xi}\rangle_{(Y_n^*,Y_n)}
-\intO{((\vr_n [\vw_n]_\delta - 2 \kappa  \nabla \mu(\vr_n)) \otimes \vv_n)(t):\Grad \vcg{\xi}}\\
&\quad+2\kappa\intO{\mu(\vr_n)\Grad\vv_n(t):\Grad \vcg{\xi}}
+ 2\kappa\intO{(\mu'(\vr_n)\vr_n - \mu(\vr_n))\Div \vv_n(t) \Div\vcg{\xi}}\\
&\quad-2 \intO{\mu(\vr_n)\Grad^t\vw_n(t):\Grad\vcg{\xi}}=0,
}
is satisfied for $\vcg{\xi}\in Y_n$, $t\in[0,T]$.
\end{itemize}
However, so far we only know that the approximate continuity equation \eqref{A_cont} is satisfied in the sense of distributions 
\eq{\label{weak1}
\intT{\langle\pt\vr_n,\phi\rangle_{(W^{-1,2}(\Omega), W^{1,2}(\Omega))}}-\intTO{\vr_n[\vw_n]_\delta\cdot\Grad\phi}+ 2\kappa\intTO{\Grad\mu(\vr_n)\cdot\Grad\phi}=0}
for any test function $\phi$ from $L^2(0,T;W^{1,2}(\Omega))$.
But on the other hand, we know that that in the distributional sense
$$\psi=\pt\vr_n- 2\kappa\Div\lr{\mu'(\vr_n)\Grad\vr_n}\in L^2((0,T),\times\Omega),$$
if only $\vw_n\in L^\infty(0,T; L^\infty(\Omega)).$  But this is the case at the level of Galerkin approximations, therefore also
$$\psi\mu'(\vr_n)\in L^2((0,T)\times\Omega).$$
Taking the product of $\psi$ and $\psi\mu'(\vr_n)$ we obtain
\eqh{\intO{\mu'(\vr_n)\lr{\pt\vr_n- 2\kappa\lap\mu(\vr_n)}^2}\leq c}
and the above integral gives rise to estimates
\eq{\intO{\mu'(\vr_n)(\pt\vr_n)^2}+ 4 \kappa^2\intO{\mu'(\vr_n)(\lap\mu(\vr_n))^2}-4\kappa\intO{\mu'(\vr_n)\pt\vr_n\lap\mu(\vr_n)}\\
=\intO{\mu'(\vr_n)(\pt\vr_n)^2}+4 \kappa^2\intO{\mu'(\vr_n)(\lap\mu(\vr_n))^2}
+2\kappa\Dt\intO{|\Grad\mu(\vr_n)|^2}\leq c.\label{h_reg}
}
Note that this estimate asks for $L^\infty((0,T)\times\Omega)$ bound for $\vw_n$, which is possible only at the level of Galerkin approximation.
 Nevertheless, regularity \eqref{h_reg} allows us to repeat (\ref{mu_w}-\ref{est}) to get the uniform (with respect to $n,\ \ep$ and $\delta$) estimate
\eq{\|\Grad\vr_n\|_{L^\infty(0,T;L^{2}(\Omega))}
+\|\lap\vr_n\|_{L^2((0,T)\times\Omega)}\leq c.
\label{high}}
\bigskip

\bigskip

\noindent{\bf Passage to the limit $n\to\infty$.} 
The biggest problem is again to pass to the limit in the term
\eq{\mu'(\vr_n)\Grad\vr_n\otimes\vw_n,\label{conv_w}}
which requires the strong convergence of the density and at least weak convergence of the gradient of density and in the convective term
\eq{\vr_n\vw_n\otimes\vw_n\label{conv_w2}}
which requires strong convergence of $\sqrt{\vr_n}\vw_n$.
Having obtained estimate \eqref{high} we can estimate the time-derivative of gradient of $\vr_n$. Indeed differentiating  \eqref{A_cont} with respect to $x$ we obtain
\eqh{\pt\Grad\vr_n=-\Grad\lr{[\vw_n]_\delta\cdot\Grad\vr_n}- 2\kappa\Grad\lap\mu(\vr_n)\in L^2(0,T; W^{-1,3/2}(\Omega)).}
Note that the above estimate is uniform also with respect to $\ep$.
Now, applying the Aubin-Lions lemma for $\Grad\vr_n$ we obtain
\eqh{\Grad\vr_n\to\Grad\vr\quad\text{strongly in } L^2(0,T;L^2(\Omega)),}
therefore due to \eqref{max} we also have
\eq{\vr_n\to\vr\quad\text{and}\quad
\frac{1}{\vr_n}\to\frac{1}{\vr}\quad\text{strongly in } L^p(0,T;L^p(\Omega))
\label{s_rho}}
for $p<\infty$ and
\eq{\vr_n\vw_n\to\vr\vw\quad \text{weakly in } L^{p_1}(0,T; L^{q_1}(\Omega))\cap L^{p_2}(0,T; L^{q_2}(\Omega)),\label{k1}}
where $ p_1<2, q_1<6, p_2<\infty,q_2<2.$
These convergences justify the limit passage in \eqref{conv_w}.

To justify passage to the limit in \eqref{conv_w2} we first  estimate 
\eq{\|\Grad\lr{\vr_n\vw_n}\|_{L^2(0,T; L^{\frac{3}{2}}(\Omega))}
\leq &\|\Grad\vr_n\|_{L^\infty(0,T;L^2(\Omega))}\|\vw_n\|_{L^2(0,T; L^6(\Omega))}\\
&+\|\Grad\vw_n\|_{L^2(0,T;L^2(\Omega))}\|\vr_n\|_{L^\infty(0,T;L^\infty(\Omega))}.
\label{k2}}
We can also estimate the time derivative of momentum, from \eqref{FG_a} we obtain
\eq{
&\sup_{\|\vcg{\phi}\|\leq1}\left|\intTO{\pt\lr{\vr_n\vw_n}\cdot\vcg{\phi}}\right|\\
&=\sup_{\|\vcg{\phi}\|\leq 1}\left\{\left|\intTO{((\vr_n [\vw_n]_\delta - 2 \kappa\mu'({\vr_n}) \nabla \vr_n) \otimes \vw_n):\Grad \vcg{\phi}}\right|\right.\\
&\quad\qquad+ 2(1-\kappa)\left|\intTO{\mu(\vr_n)D(\vw_n):\Grad \vcg{\phi}}\right|\\
&\quad\qquad+\ep\left|\intTOB{\lap\vw\cdot\lap\vcg{\phi}}\right|
+\ep\left|\intTOB{(1+|\Grad\vw|^2) \Grad\vw:\Grad\vcg{\phi}}\right|
\\
&\quad\qquad+\left.
2 \kappa\left|\intTO{\mu(\vr_n)A(\vw_n): \Grad\vcg{\phi}}\right|
+2\kappa(1-\kappa)\left|\intTO{\mu(\vr_n)\Grad\vv_n:\Grad\vcg{\phi}}\right|\right\},
\label{ptw}}
where $\|\vcg{\phi}\|$ denotes the norm in the space $W_T:=L^2(0,T;V\cap W^{2,2}(\Omega))\cap L^4(0,T;W^{1,4}(\Omega))$.
Let us estimate the convective term
\eqh{&\left|\intTO{((\vr_n [\vw_n]_\delta - 2 \kappa\mu'({\vr_n}) \nabla \vr_n) \otimes \vw_n):\Grad \vcg{\phi}}\right|\\
&\leq\intT{\|\Grad\vcg{\phi}\|_{L^6(\Omega)}\lr{R\|\vw_n\|^2_{L^\frac{12}{5}(\Omega)}+c(\kappa,R)\|\Grad\vr_n\|_{L^6(\Omega)}\|\vw_n\|_{L^\frac{3}{2}(\Omega)}}}\\
&\leq c(\kappa,R,\ep)\|\vcg{\phi}\|_{L^2(0,T; W^{2,2}(\Omega))},
}
for the highest order terms we have
\eqh{\ep\left|\intTO{\lap\vw_n\cdot\lap\vcg{\phi}}\right|\leq \ep\|\vw_n\|_{L^2(0,T; W^{2,2}(\Omega))}\|\vcg{\phi}\|_{L^2(0,T; W^{2,2}(\Omega)),}}
and
\eqh{\ep\left|\intTO{(1+|\Grad\vw|^2) \Grad\vw:\Grad\vcg{\phi}}\right|
\leq \ep\intT{\|\Grad\vcg{\phi}\|_{L^4(\Omega)}\lr{\|\Grad\vw_n\|^3_{L^{4}(\Omega)}+\|\Grad\vw_n\|_{L^{\frac{4}{3}}(\Omega)}}}\\
\leq \ep\|\Grad\vcg{\phi}\|_{L^4(0,T; L^4(\Omega))}\lr{\|\Grad\vw\|^3_{L^4(0,T; L^4(\Omega))}+\|\Grad\vw\|_{L^{\frac{4}{3}}(0,T; L^{\frac{4}{3}}(\Omega))}}
.}
The other terms from \eqref{ptw} are less restrictive, therefore
\eq{\|\pt\lr{\vr_n\vw_n}\|_{W_T^*}\leq c,\label{k3}}
where $W_T^*$ denotes the dual space of $W_T$ defined above.
Collecting \eqref{k1}, \eqref{k2}, \eqref{k3} and applying the Aubin-Lions lemma to $\vr_n\vw_n$, we prove that
\eqh{\vr_n\vw_n\to\vr\vw\quad \text{strongly in } L^p(0,T;L^p(\Omega))}
for some $p>1$ and therefore thanks to \eqref{s_rho} and \eqref{ueu}
\eq{\Grad\vw_n\to\Grad\vw\quad \text{strongly in }L^p(0,T;L^p(\Omega))\label{strong_gw}}
 for $1\leq p<4$. In particular, convergence in \eqref{conv_w2} is  proved. 
     For future purposes we now estimate the time derivative of $\vr\vv$ in $L^{\frac{4}{3}}(0,T; W^{-1,\frac{4}{3}}(\Omega))$. We use \eqref{FGtheta_a}
to obtain
 \eqh{
&\sup_{\|\vcg{\xi}\|\leq1}\left|\intTO{\pt\lr{\vr_n\vv_n}\cdot{\vcg{\xi}}}\right|\\
&= \sup_{\|\vcg{\xi}\|\leq1}\left\{\left|\intTO{((\vr_n [\vw_n]_\delta - 2 \kappa \mu'({\vr_n}) \nabla \vr_n) \otimes \vv_n):\Grad \vcg{\xi}}\right|\right.\\
&\qquad+ 2 \kappa\left|\intTO{\mu(\vr_n)\Grad\vv_n:\Grad \vcg{\xi}}\right|
+ 2\kappa\left|\intTO{(\mu'(\vr_n)\vr_n - \mu(\vr_n))\Div \vv_n\ \Div\vcg{\xi}}\right|\\
&\qquad+ 2\left.\left|\intTO{\mu(\vr_n)\Grad^t\vw_n:\Grad\vcg{\xi}}\right|\right\}}
 where $\|\vcg{\xi}\|$ denotes the norm in the space $L^4(0,T;W^{1,4}(\Omega))$. 
 We will only estimate the convective term since it is most restrictive.
\eqh{ 
&\left|\intTO{((\vr_n [\vw_n]_\delta - 2 \kappa \mu'({\vr_n}) \nabla \vr_n) \otimes \vv_n):\Grad \vcg{\xi}}\right|\\
&\leq\intT{\|\Grad\vcg{\xi}\|_{L^4(\Omega)}\lr{R\|\vw_n\|_{L^4(\Omega)}\|\vv_n\|_{L^2(\Omega)}+c(\kappa,R)\|\Grad\vr_n\|_{L^4(\Omega)}\|\vv_n\|_{L^2(\Omega)}}}\\
&\leq c(\kappa,R)\|\vcg{\xi}\|_{L^4(0,T; W^{1,4}(\Omega))}\|\vv_n\|_{L^2(0,T;L^2(\Omega))}
\lr{\|\vw_n\|_{L^4(0,T; L^4(\Omega))}+\|\lap\vr_n\|^{\frac{1}{2}}_{L^2(0,T; L^2(\Omega))}}
}
thus
\eq{\|\pt\lr{\vr_n\vv_n}\|_{L^{\frac{4}{3}}(0,T; W^{-1,\frac{4}{3}}(\Omega))}\leq c.\label{k4}}
Hence, the limit functions $(\vr,\vw,\vv)=(\vr_\delta,\vw_\delta,\vv_\delta)$ fulfil
\begin{itemize}
\item the continuity equation 
\eq{\pt\vr_\delta+\Div\lr{\vr_\delta[\vw_\delta]_\delta}- 2\kappa\lap\mu(\vr_\delta)=0\label{strong1}}
 a.e. in $(0,T)\times\Omega$,
\item the momentum equation
\begin{equation}\label{FGn}
\begin{split}
&
\langle\pt\lr{\vr_\delta\vw_\delta},\vcg{\phi}\rangle_{(W_\tau^*, W_\tau)}
-\inttauO{((\vr_\delta [\vw_\delta]_\delta - 2 \kappa \nabla \mu(\vr_\delta)) \otimes \vw_\delta):\Grad \vcg{\phi}}\\
&\quad+ 2(1-\kappa)\inttauO{\mu(\vr_\delta)D(\vw_\delta):\Grad \vcg{\phi}}\\
&\quad+\ep\inttauO{\lr{\lap\vw_\delta\cdot\lap\vcg{\phi}}}
+\ep\inttauO{\lr{(1+|\Grad\vw_\delta|^2) \Grad\vw_\delta:\Grad\vcg{\phi}}}\\
&\quad+2 \kappa\inttauO{\mu(\vr_\delta)A(\vw_\delta): \Grad\vcg{\phi}} -2 \kappa(1-\kappa)\inttauO{\mu(\vr_\delta)\Grad\vv_\delta:\Grad\vcg{\phi}}=0,
\end{split}
\end{equation}
for all $\vcg{\phi}\in W_\tau$ with $\tau\in[0,T]$,
\item the auxiliary equation for $\vv$
\eq{\label{FGthetan}
&
\langle\pt\lr{\vr_\delta\vv_\delta},\vcg{\xi}\rangle_{(L^{\frac{4}{3}}(0,\tau; W^{-1,\frac{4}{3}}(\Omega)), L^4(0,\tau; W^{1,4}(\Omega)))}
\\
&\quad-\inttauO{((\vr_\delta [\vw_\delta]_\delta - 2 \kappa  \nabla \mu(\vr_\delta)) \otimes \vv_\delta):\Grad \vcg{\xi}}+ 2\kappa\inttauO{\mu(\vr_\delta)\Grad\vv_\delta:\Grad \vcg{\xi}}\\
&\quad+ 2\kappa\inttauO{(\mu'(\vr_\delta)\vr_\delta - \mu(\vr_\delta))\Div \vv_\delta\ \Div\vcg{\xi}}-2 \inttauO{\mu(\vr_\delta)\Grad^t\vw_\delta:\Grad\vcg{\xi}}=0,
}
for all $\vcg{\xi}\in L^4(0,\tau; W^{1,4}(\Omega))$ with $\tau\in[0,T]$.
\end{itemize}

\bigskip

\noindent{\bf Passage to the limit $\delta$ tends to $0$ and identification of $\vv_\delta$ with $\Grad\vp(\vr_\delta)$ at the limit.}
The aim of this paragraph is to let $\delta\to0$ in the equations \eqref{strong1}, \eqref{FGn} and  \eqref{FGn}. This limit passage can be performed exactly as $n\to \infty$ presented above. The only difference is that after this step we may drop the additional equation for $\vv$ thanks to identification $\vv=2 \Grad\vp(\vr)$. Below we present the details of this reasoning.

Note that the coefficients of the quasi-linear parabolic equation \eqref{strong1} (i.e. $[\vw_\delta]_\delta$) are sufficiently regular and the maximum principle \eqref{max} holds uniformly with respect to all approximation parameters. Therefore, the classical theory  of Lady{\v{z}}enskaja, Solonnikov and Uralceva \cite{LSU} (Theorems 7.2, 7.3 and 7.4 from \cite{LSU}) can be applied to show further regularity of $\vr_\delta$, we have in particular
\eq{\pt\vr_\delta\in C([0,T]; C(\Omega)),\quad \vr_\delta\in C([0,T]; C^2(\Omega)).\label{reg_rg}}
Let us now rewrite \eqref{strong1} using \eqref{rel1} as
$$  \pt\vr_\delta
     + \Div\lr{\vr_\delta  ([\vw_\delta]_\delta - 2 \kappa\, \nabla\vp(\vr_\delta)) } = 0$$  
and therefore multiplying the above equation by $\mu'(\vr)$ we obtain 
$$\partial_t\mu(\vr_\delta) 
     + \Div\lr{\mu(\vr_\delta)([\vw_\delta]_\delta 
       - 2 \kappa\, \Grad\vp(\vr_\delta))} 
     - 2\kappa (\mu'(\vr_\delta)\vr_\delta - \mu(\vr_\delta)) \Delta \vp(\vr_\delta)  
     = 0. $$  
Differentiating it with respect to space one gets in the sense of distributions
\eq{ \pt(\vr_\delta \widetilde \vv_\delta) 
     + \Div(\vr_\delta ([\vw_\delta]_\delta 
     -2 \kappa\, \Grad\vp(\vr_\delta))\otimes \widetilde\vv_\delta ) 
     -2\kappa\Grad\lr{ (\mu'(\vr_\delta)\vr_\delta - \mu(\vr_\delta)) 
        \Div \widetilde\vv_\delta  }\\
    +2 \Div((\mu(\vr_\delta) \Grad^{t} [\vw_\delta]_\delta)
    - 2 \kappa \Div(\mu(\vr_\delta)\Grad\widetilde \vv_\delta)
      = 0\label{test}}
 where by $\widetilde\vv_\delta$ we denoted $2\Grad\vp(\vr_\delta)$. 
 Note that due to particular case of the Gagliardo-Nirenberg interpolation inequality \eqref{G-N} and \eqref{max} we know that $\Grad\vr_\delta$ is bounded
  \eqh{\|\nabla \vr_\delta\|_{L^4(0,T; L^4(\Omega))}\leq c,} 
 uniformly with respect to $\delta$.  One can thus estimate the convective term of \eqref{test} in $L^2(0,T; W^{-1,2}(\Omega))$ uniformly with respect to $\delta$. Indeed, we now $\vw_\delta$ uniformly
bounded in $L^4(0,T;L^4(\Omega))$ with respect to $\delta$ and therefore
\eqh{&\sup_{\|\vcg{\xi}\|\leq1}\left|\intTO{  (\vr_\delta ([\vw_\delta]_\delta 
     - 2 \kappa\, \Grad\vp(\vr_\delta))\otimes \widetilde\vv_\delta ):\Grad\vcg{\xi}}\right|\\
     &\leq c(R) \|\Grad\vcg{\xi}\|_{L^2(0,T; L^2(\Omega))}\|\nabla \vr_\delta\|_{L^4(0,T; L^4(\Omega))}
     \lr{\|\nabla [\vw_\delta]_\delta\|_{L^4(0,T; L^4(\Omega))}+\|\nabla \vr_\delta\|_{L^4(0,T; L^4(\Omega))}}
     }
 for $\vcg{\xi}\in L^2(0,T; W^{1,2}(\Omega))$ (uniformly with respect to $\delta$), which justifies that 
\eq{ &\langle \pt(\vr_\delta \widetilde \vv_\delta), \vcg{\xi}\rangle_{(L^2(0,T; W^{-1,2}(\Omega)),L^2(0,T; W^{1,2}(\Omega)))}\\
     &-\intTO{\vr_\delta ([\vw_\delta]_\delta 
     - 2 \kappa\, \Grad\vp(\vr_\delta))\otimes \widetilde\vv_\delta :\Grad\vcg{\xi}} 
     + 2\kappa\intTO{ (\mu'(\vr_\delta)\vr_\delta - \mu(\vr_\delta))  \Div \widetilde\vv_\delta \Div\vcg{\xi} }\\
   &  - 2\intTO{\mu(\vr_\delta) \Grad^{t} [\vw_\delta]_\delta:\Grad\vcg{\xi}}
    + 2\kappa \intTO{\mu(\vr_\delta)\Grad\widetilde \vv_\delta:\Grad\vcg{\xi}}
      = 0\label{test2}}
 is satisfied for any $\vcg{\xi}\in L^2(0,T; W^{1,2}(\Omega))$.
 \begin{rmk}
 As we noticed in \eqref{reg_rg} the regularity of $\vr_\delta$ is in fact much higher and allows to formulate the equation for $\widetilde\vv_\delta$ in much stronger sense than merely \eqref{test2}. This formulation, however, will be used when passing to the limit with respect to $\ep$ after passing to the limit with respect to $\delta$.
 \end{rmk}
 We now want  show that  $\widetilde\vv_\delta -\vv_\delta$ tends to $0$ when $\delta$ goes to zero
 in an appropriate norm. To this purpose let us expand
  \eq{
I= & \Dt\intO{\vr_\delta\frac{|\vv_\delta-\widetilde\vv_\delta|^2}{2}}
  + 2 \kappa \intO{\mu(\vr_\delta)|\nabla (\vv_\delta - \widetilde\vv_\delta)|^2} \\
  & + 2\kappa\intTO{(\mu'(\vr_\delta)\vr_\delta - \mu(\vr_\delta))(\Div(\vv_\delta -\widetilde \vv_\delta))^2}\\
 =&  \Dt\intO{\vr\lr{\frac{|\vv_\delta|^2}{2}-\vv_\delta \cdot\widetilde\vv_\delta+\frac{|\widetilde\vv_\delta|^2}{2}}}\\ 
 &+ 2\kappa\intO{\mu(\vr_\delta)\lr{|\nabla\vv_\delta |^2 +|\nabla \widetilde\vv_\delta|^2
    - 2 \nabla \vv_\delta\cdot \nabla \widetilde \vv_\delta}}. \\
 & +2 \kappa\intTO{(\mu'(\vr_\delta)\vr_\delta - \mu(\vr_\delta))\lr{|\Div \vv_\delta|^2 + |\Div \widetilde \vv_\delta|^2  - 2 \Div \vv_\delta \Div \widetilde \vv_\delta}}
  \label{expand}
 }
 To handle the first term, let us notice that letting $n\to \infty$ in \eqref{t1}, using the lower semi-continuity of the convex functions and the strong convergence of $\Grad\vw_n$ established in \eqref{strong_gw} we obtain
 \eq{\label{t11}
&\Dt\intO{\vr_\delta\frac{|\vv_\delta|^2}{2}}
+ 2\kappa\intTO{\mu(\vr_\delta)|\Grad\vv_\delta|^2}\\
&\quad+ 2\kappa\intTO{(\mu'(\vr_\delta)\vr_\delta - \mu(\vr_\delta))(\Div \vv_\delta)^2}-
2 \intTO{\mu(\vr_\delta)\Grad^t\vw_\delta:\Grad\vv_\delta}\leq0.
}
 Now, the last term in \eqref{expand} can be computed using $\vcg{\xi}=\widetilde\vv_\delta$ in \eqref{test2}, we have
  \eqh{&\Dt\intO{\vr_\delta\frac{|\widetilde\vv_\delta|^2}{2}} 
   +2 \kappa\intO{ (\mu'(\vr_\delta)\vr_\delta - \mu(\vr_\delta))  (\Div \widetilde\vv_\delta)^2 }\\
   & - 2\intO{\mu(\vr_\delta) \Grad^{t} [\vw_\delta]_\delta :\Grad\widetilde\vv_\delta}
    + 2\kappa \intO{\mu(\vr_\delta)|\Grad\widetilde \vv_\delta|^2}=0.
 }
The middle term in \eqref{expand} equals
\eq{\Dt\intO{\vr_\delta\vv_\delta\cdot\widetilde\vv_\delta}=\intOB{\pt\lr{\vr_\delta\vv_\delta}\cdot\widetilde\vv_\delta + \vv_\delta \cdot\pt(\vr_\delta \widetilde\vv_\delta ) - \partial_t\vr_\delta\,  \vv_\delta\cdot \widetilde \vv_\delta}\label{der_mixed}}
and the two first terms make sense and can be handled using $\vcg{\xi}=\widetilde\vv_\delta$ in \eqref{FGthetan} and $\vcg{\xi} = \vv_\delta$ in \eqref{test2}. 
Note that $\widetilde \vv_\delta$ and $\partial_t \vr_\delta$ are due to \eqref{reg_rg} regular enough to justify the integrability of the last term in \eqref{der_mixed} and we can write
\eqh{&\intO{\pt\vr_\delta\vv_\delta\cdot\widetilde\vv_\delta}\\
&=\intO{\lr{\vr_\delta  ([\vw_\delta]_\delta - 2 \kappa\, \nabla\vp(\vr_\delta)) } \otimes\vv_\delta:\Grad\widetilde\vv_\delta}
+\intO{\lr{\vr_\delta  ([\vw_\delta]_\delta - 2 \kappa\, \nabla\vp(\vr_\delta)) } \otimes\widetilde\vv_\delta:\Grad\vv_\delta}.
} 
Therefore, after summing all expressions together and after some manipulation, we can show that
\eqh{I-\intO{\mu(\vr_\delta)\Grad^t[\vw_\delta]_\delta:\Grad\widetilde\vv_\delta}
-\intO{\mu(\vr_\delta)\Grad^t\vw_\delta:\Grad\vv_\delta}\\
+\intO{\mu(\vr_\delta)\Grad^t\vw_\delta:\Grad\widetilde\vv_\delta}
+\intO{\mu(\vr_\delta)\Grad^t[\vw_\delta]_\delta:\Grad\vv_\delta}\leq 0,
}
in particular
$$I \leq \intO{\mu(\vr)\Grad^t([\vw_\delta]_\delta - \vw_\delta ): \Grad (\widetilde\vv_\delta - \vv_\delta)}.$$
Note that the r.h.s. of this inequality tends to 0 when $\delta\to 0$. Indeed, we can bound $ \Grad (\widetilde\vv_\delta - \vv_\delta)$ in $L^2(0,T; L^2(\Omega))$ uniformly with respect to $\delta$ and $[\vw_\delta]_\delta\to \vw$ strongly in $L^p(0,T; W^{1,p}(\Omega))$ for $p<4$.
Therefore, using \eqref{expand}, we conclude that $\vv_\delta - \widetilde\vv_\delta$ converges to zero in $L^\infty(0,T;L^2(\Omega))\cap L^2(0,T;H^1(\Omega))$ when $\delta \to 0$. $\Box$

The limit functions $(\vr,\vw)=(\vr_\ep,\vw_\ep)$ fulfil
\begin{itemize}
\item the continuity equation 
\eq{\pt\vr_\ep+\Div\lr{\vr_\ep\vw_\ep}- 2\kappa\lap\mu(\vr_\ep)=0\label{cont_ep}}
 a.e. in $(0,T)\times\Omega$,
\item the momentum equation
\begin{equation}\label{FG_ep}
\begin{split}
&
\langle\pt\lr{\vr_\ep\vw_\ep},\vcg{\phi}\rangle_{(W_\tau^*, W_\tau)}
-\inttauO{((\vr_\ep\vw_\ep -  2\kappa \nabla \mu(\vr_\ep)) \otimes \vw_\ep):\Grad \vcg{\phi}}\\
&\quad+ 2(1-\kappa)\inttauO{\mu(\vr_\ep)D(\vw_\delta):\Grad \vcg{\phi}}\\
&\quad+\ep\inttauO{\lr{\lap\vw_\ep\cdot\lap\vcg{\phi}}}
+\ep\inttauO{\lr{(1+|\Grad\vw_\ep|^2) \Grad\vw_\ep:\Grad\vcg{\phi}}}\\
&\quad+2 \kappa\inttauO{\mu(\vr_\ep)A(\vw_\ep): \Grad\vcg{\phi}}- 4 \kappa(1-\kappa)\inttauO{\mu(\vr_\ep)\Grad^2\vp(\vr_\ep):\Grad\vcg{\phi}}=0,
\end{split}
\end{equation}
for $\vcg{\phi}\in W_\tau$, $\tau\in[0,T]$,
\item the auxiliary equation for $\Grad\vp(\vr_\ep)$
\eq{\label{FGtheta_ep}
&
\langle\pt\Grad\mu(\vr_\ep),\vcg{\xi}\rangle_{(L^2(0,\tau; W^{-1,2}(\Omega)),L^2(0,\tau; W^{1,2}(\Omega)))}
\\
&\quad-\inttauO{((\vr_\ep\vw_\ep-  2\kappa  \nabla \mu(\vr_\ep)) \otimes \Grad\vp(\vr_\ep)):\Grad \vcg{\xi}}+2\kappa\inttauO{\Grad\mu(\vr_\ep):\Grad \vcg{\xi}}\\
&\quad+2\kappa\inttauO{(\mu'(\vr_\ep)\vr_\ep - \mu(\vr_\ep))\lap\vp(\vr_\ep)\ \Div\vcg{\xi}}-\inttauO{\mu(\vr_\ep)\Grad^t\vw_\ep:\Grad\vcg{\xi}}=0,
}
for $\vcg{\xi}\in L^2(0,\tau; W^{1,2}(\Omega))$, $\tau\in[0,T]$.
\end{itemize}

\bigskip

\noindent {\bf Passage to the limit with respect to $\varepsilon$: existence of global solutions.}
Let us start this paragraph by recalling the estimates that are uniform with respect to $\ep$.
Passing to the limit $\delta\to 0$ in  \eqref{aa1}, using the weak convergence of $\Grad\vr_\delta$ and strong convergence of $\vr_\ep$ together with strong convergence of $\Grad\vw_\ep$ and weak convergence of $\lap\vw_\ep$ and a standard argument based on convexity of norm we obtain
\eq{\label{aa2}
\Dt&\intO{\vr_\ep\lr{\frac{|\vw_\ep|^2}{2}
+(1-\kappa)\kappa\frac{|2\nabla\vp(\vr_\ep)|^2}{2}}}+2(1-\kappa)\intO{\mu(\vr) |D(\vw_\ep) -2\kappa\nabla \Grad\vp(\vr_\ep))|^2} \\
&+\ep\intOB{|\lap\vw_\ep|^2+(1+|\Grad\vw_\ep|^2) |\Grad\vw_\ep|^2}\\
&+ 2\kappa\intO{\mu(\vr_\ep) |A(\vw_\ep)|^2}+ 2(1-\kappa)\intO{(\mu'(\vr_\ep)\vr_\ep-\mu(\vr_\ep))|2\kappa \lap\vp(\vr_\ep)|^2}\leq0.
}
We also have \eqref{max} and thus we can deduce from the above inequality that
\eq{\|\vw_\ep\|_{L^\infty(0,T;H)}+\|\vw_\ep\|_{L^2(0,T;V)}
+\ep^{1/2}\|\vw_\ep\|_{L^2(0,T;V\cap W^{2,2}(\Omega))}+
\ep^{1/4}\|\Grad\vw_\ep\|_{L^4(0,T;L^4(\Omega))}\\
+\|\Grad\vr_\ep\|_{L^\infty(0,T;L^2(\Omega))}
+\|\vr_\ep\|_{L^2(0,T;W^{2,2}(\Omega))}\leq c.
\label{v3}}

 Having obtained estimate \eqref{v3} we are ready to perform the last limit passage and to deduce existence of weak solutions to original system \eqref{main1} in the sense of Definition \eqref{Def1}.  This in fact can be done almost exactly the same way as the proof of sequential stability of solutions from the paper \cite{BrEsSy07}. The only difference is to pass to the limit in the term
$\Grad\mu(\vr)\otimes\Grad\vp(\vr)$ in \eqref{FG_ep} and \eqref{FGtheta_ep}. To this purpose we need the following compensated compactness lemma (see for instance \cite{PLL} Lemma 5.1).
\begin{lemma}\label{lem_compac} Let $g_n$, $h_n$ converge weakly to $g$, $h$ respectively in $L^{p_1}(0,T;L^{p_2}(\Omega))$, $L^{q_1}(0,T;L^{q_2}(\Omega))$ where $1\leq p_1,p_2 \leq \infty$, $\frac{1}{p_1}+ \frac{1}{q_1}=\frac{1}{p_2}+ \frac{1}{q_2}= 1$. We assume in addition that
\begin{itemize}
  \item $\partial_t g_n$ is bounded in $L^1(0,T;W^{-m,1}(\Omega))$ for some $m\geq 0$ independent of $n$.
  \item $\|h_n-h_n(t,\cdot+\xi)\|_{L^{q_1}(L^{q_2})} \rightarrow 0$ as $|\xi|\rightarrow 0$, uniformly in $n$.
\end{itemize}
Then $g_nh_n$ converges to $gh$ in $\mathcal{D}'$.
\end{lemma}
   We will apply this lemma to $g_n=\Grad\mu(\vr_\ep)$ and $h_n=\Grad\vp(\vr_\ep)$. First of all, let us notice that exactly as in the previous sections we have $\vr_\ep\to\vr$ strongly in $L^p((0,T)\times\Omega)$ for any $p$ finite, and so we have also that $\mu(\vr_\ep), \vp(\vr_\ep)$ converge strongly to $\mu(\vr),\ \vp(\vr)$, respectively. Due to \eqref{v3} both $g_n$ and $h_n$ are bounded in $L^2((0,T)\times\Omega)$ and they converge weakly to $g=\Grad\mu(\vr)$, $h=\Grad\vp(\vr)$. From the same estimate it follows that $\Grad h_n$ is uniformly bounded in $L^2((0,T)\times\Omega)$. 

Moreover, the convective term in \eqref{FGtheta_ep} is bounded uniformly with respect to $\ep$ for $\vcg{\xi}\in L^2(0,T; W^{1,3}(\Omega))$ therefore we may estimate $\pt g_n$ in $L^2(0,T; W^{-1,\frac{3}{2}}(\Omega))$, thus the product $g_n h_n$ converges to $\Grad\mu(\vr)\otimes\Grad\vp(\vr)$ in the sense of distributions on $(0,T)\times\Omega$.
   
  Finally, the passage $\ep\to0$ in all $\ep$-dependent terms of \eqref{FG_ep}  gives $0$ due to the uniform estimates from \eqref{v3}. 
  To conclude let us check in which sense are the initial conditions admitted, Since the time derivative of $\Grad\vr_\ep$ is bounded in $L^2(0,T; W^{-1,\frac{3}{2}}(\Omega))$ and $\Grad\vr\in L^\infty(0,T; L^2(\Omega))$, we may use the Arzela-Ascoli\`e theorem to verify that $\Grad\vr_\ep\to\Grad\vr$ in  $C([0,T]; L^2_{\rm weak}(\Omega))$ \footnote{$\zeta\in C([0,T]; L^2_{\rm weak}(\Omega))$ iff $\lim_{t\to t_0}|\langle{\eta,\zeta(t)-\zeta(t_0)\rangle}|=0$, $\forall \eta\in L^2(\Omega)$, $\forall t_0\in[0,T]$}, also $\vr_\ep\vw_\ep\to\vr\vw$ in $C([0,T]; L^2_{\rm weak}(\Omega))$. Moreover, using a version of Aubin-Lions lemma we obtain $\vr$ is strongly continuous, i.e. $\vr_\ep\to\vr$ in $C([0,T]; L^2(\Omega))$. $\Box$

.

\section{Proof of Theorem \ref{T_main2}}
\subsection{Part 1, passage to the limit $\kappa\to 0$}
The aim of this section will be to let $\kappa\to 0$ in \eqref{main1} and to show that the sequence of solutions $(\vr_\kappa,\vw_\kappa)$ converges to $(\vr,\vu)$ the weak solution of  the non-homogeneous incompressible Navier-Stokes system \eqref{mainINS}
\begin{equation}\label{main2}
\begin{array}{c}
\pt\vr+ \Div(\vr \vu)  = 0,\\
\pt\lr{\vr \vu}+ \Div  (\vr \vu \otimes \vu) - 2 \Div (\mu(\vr) D(\vu)) + \Grad {\pi_1} =\vc{0},\\
\Div\vu = 0.
\end{array}
\end{equation}

\noindent {\bf Strong convergence of the density.} We start the proof of Theorem \ref{T_main2} by proving that the uniform estimates obtained in the previous section can be used to deduce the strong convergence of the density. First let us recall that uniformly with respect to $\kappa$ we have
\eq{0<r\leq\vr_\kappa \leq R<\infty.\label{max_vr_k}}
For any $\kappa$ fixed we have
\eqh{\vw_\kappa\in{L^2(0,T; V)}\cap L^\infty(0,T;H).}
However, the bounds from \eqref{v2} allow us only to show that
\eqh{
&\vu_\kappa\to\vu\quad  \text{weakly\ in\ } L^2(0,T;H^1(\Omega))\\
&   \kappa\vr_\kappa\to 0 \quad \text{strongly\ in\ }L^\infty(0,T;H^1(\Omega))\cap L^2(0,T; H^2(\Omega)),\\
&\vw_\kappa\to\vw\quad \text{weakly$^*$\ in\ } L^\infty(0,T;H).\label{weak_w}
}
Thus it follows that
\eq{
&   \vw_\kappa=\vu_\kappa+ 2\kappa \Grad \vp(\vr_\kappa)\to \vu \quad \text{weakly\ in\ }L^2(0,T;H^1(\Omega)),\\
&\vw_\kappa\to\vu\quad \text{weakly$^*$\ in\ } L^\infty(0,T;H),\label{weak_w}}
in particular
\eq{\vw_\kappa\to \vu \quad \text{weakly\ in\ }L^2(0,T;V).}
Moreover, we know that the pair $(\vr_\kappa,\vw_\kappa)$ satisfies the equation
\eq{\pt \vr_\kappa + \vw_\kappa\cdot\Grad \vr_\kappa - 2\kappa \lap\mu(\vr_\kappa) = 0\label{ca}}
in the sense of distributions on $(0,T)\times\Omega$.
As a consequence, one may deduce using Aubin-Lions lemma and Arzela-Ascoli\`e theorem that there exists a subsequence s.t.
\eq{\label{conv}
\vr_\kappa\to\vr \quad \text{in\ } C([0,T]; W^{-m,p}(\Omega)),\quad\text{and}\quad \vr_\kappa\to\vr \quad \text{in\ } C([0,T]; L^{p}_{\rm{weak}}(\Omega))}
for $m>0$, $1\leq p<\infty$. Using this and \eqref{weak_w} we prove that
\eqh{\vr_\kappa\vw_\kappa\to \vr\vu,\quad\text{weakly\ in\ } L^p(0,T;(L^p(\Omega))^3)}
for some $p>1$. Therefore the limit pair $(\vr,\vu)$ satisfies the continuity equation
\eq{\pt \vr + \vu\cdot\Grad \vr = 0\label{lim_cont}}
in the sense of distributions on $(0,T)\times\Omega$ and
$$\Div\vu=0\quad  \text{a.e. in }(0,T)\times\Omega.$$
Now, our aim will be to prove strong convergence of the density, which is given by the following lemma
\begin{lemma}\label{L_strong}
Let $(\vr_\kappa,\vw_\kappa)_{\kappa>0}$ be a sequence of solutions satisfying the above weak convergences, then $\vr_\kappa$ converges to $\vr$ strongly in $C([0,T];L^p(\Omega))$ for all $1\leq p<\infty$. Moreover, $\vr$ is the unique solution of \eqref{lim_cont}.
\end{lemma}
\pf Uniqueness of $\vr$  and the fact that $\vr$ belongs to $C([0,T];L^p(\Omega))$ relies on the fact that if $\vu\in L^2(0,T;V)$, $\vr\in L^\infty((0,T)\times\Omega)$ and \eqref{lim_cont} is satisfied, then for any $\beta\in C(\mathbb{R})$,  $\beta(\vr)$ also solves \eqref{lim_cont}.
This can proven exactly as in P.-L.~Lions \cite{PLL} Theorem 2.4., Step 2, thus the task is only to show the strong convergence of the density.
For this purpose, let us observe that since $\vr,\vu$ satisfy \eqref{lim_cont}, taking $\beta(\vr)=\vr^2$ we obtain the equation
\eq{\pt \vr^2 + \vu\cdot\Grad\vr^2 = 0\label{lim_cont_ren},}
therefore integrating over $(0,T)\times\Omega$ we check that
\eq{\label{vrlim}
\intO{\lr{\vr(t)}^2}=\intO{\lr{\vr^0}^2}
}
for all $t\in[0,T]$.
Next, multiplying the approximate equation \eqref{ca} by $2\vr_\kappa$ and integrating by parts, we obtain
\eq{\label{vrkappa}
\intO{\lr{\vr_\kappa (t)}^2}=-4\kappa\int_0^t\intO{\mu'(\vr_k)|\Grad\vr_\kappa|^2}\, d{\rm s}+\intO{\lr{\vr_\kappa^0}^2}
}
Taking the sup in time, passing to the limit and using the strong convergence of initial data we therefore obtain
\eqh{\limsup_{\kappa\to0} \sup_{t\in [0,T]}
\intO{\lr{\vr_\kappa (t)}^2}\leq \intO{\lr{\vr^0}^2}= \sup_{t\in [0,T]} \intO{\lr{\vr(t)}^2}.
}
Note however, that by the lower semicontinuity of norms we know that
$$ \sup_{t\in [0,T]} \intO{(\vr(t))^2}  \le \liminf_{\kappa\to 0} \sup_{t\in [0,T]} \intO{(\vr_\kappa(t)) ^2}.$$
Combining these two inequalities,  we finally obtain
\eqh{\lim_{\kappa\to0} \sup_{t\in [0,T]}
\intO{\lr{\vr_\kappa (t)}^2}= \sup_{t\in [0,T]} \intO{\lr{\vr(t)}^2}.
}
Coming back to \eqref{vrkappa}, we get the strong convergence 
$\kappa \int_0^t\intO {\mu'(\vr_k) |\nabla \vr_\kappa|^2}d\tau \to 0 \hbox{ for all } t \in (0,T]$ 
and therefore coming back to \eqref{vrkappa}   
\eqh{\lim_{\kappa\to0}
\intO{\lr{\vr_\kappa (t)}^2}=  \intO{\lr{\vr(t)}^2}.
}
for all $t\in[0,T]$ on account of \eqref{vrlim}. 
 In view of the convergence in $C([0,T];L^2_{{\rm weak}}(\Omega))$, this 
implies that  
$$\vr_\kappa\to \vr\quad \text{in\ } C([0,T]; L^{2}(\Omega))$$
and the statement of Lemma \ref{L_strong} follows by \eqref{conv}. $\Box$

\bigskip

\noindent {\bf Passage to the limit in the momentum equation.} On account of \eqref{weak_w} and strong convergence of the density stablished above, the passage to the limit in \eqref{weak_mom} requires checking the limit in the nonlinear terms $\vr_\kappa\vw_\kappa\otimes\vw_\kappa$ and $\kappa\Grad\mu(\vr_\kappa)\otimes\Grad\vp(\vr_\kappa)$. The passage to the limit in the first term can be justified by application of Lemma \ref{lem_compac} with $g_n=\vr_\kappa\vw_\kappa$, $h_n=\vw_\kappa$, the details are left to the reader.
To pass to the limit in the second term we only need to show that $\kappa|\Grad\vr_\kappa|^2$ converges to $0$. Note that it does not follow from the estimates obtained in \eqref{v2} which only provide  uniform $L^1$ bound of the term in question. To solve this problem we need to recall that  we have proved previously that
\eqh{\lim_{\kappa\to0}\kappa\intTO{\mu'(\vr_\kappa)|\Grad\vr_\kappa|^2}=0}
which due to \eqref{max_vr_k} implies that 
\eqh{\kappa|\Grad\vr_\kappa|^2\to0\quad\text{strongly in}\ L^1((0,T)\times\Omega).}
Therefore, the limit momentum equation reads
\eqh{
&\intTO{\vr\vu\cdot\pt\vphi}+\intO{\vr\vu\otimes\vu:\Grad\vphi}-2\intTO{\mu(\vr)D(\vu):\Grad\vphi}=-\intO{\vu^0\cdot\vphi(0)}}
and is satisfied for $\vphi\in (C^\infty((0,T)\times\Omega))^3$, s.t.  $\Div\vphi=0$ and $\vphi(T)=\vc{0}$.

\subsection{Part 2, passage to the limit $\kappa\to 1$} The aim of this section will be to let $\kappa\to 1$ in \eqref{main1} and to show that the sequence of solutions $(\vr_\kappa,\vw_\kappa)$ converges to the global weak solution  $(\vr,\vu)$ of the Kahzikhov-Smagulov type system,  i.e. \eqref{main2}.

   The basic idea is again to use the estimate \eqref{v2}, \eqref{max_vr_k}, \eqref{est} and  deduce the convergences
\begin{equation}\label{conv2}
\begin{gathered}
\vw_\kappa\rightharpoonup \vw  \quad \text{weakly in } L^2(0,T;V) \text{ and weakly$^*$ in } 
     L^\infty(0,T;H),\\
     \vr_\kappa\to\vr\quad \text{weakly\ in\ } L^p((0,T)\times\Omega),\ 1\leq p<\infty,\text{ and weakly in } 
     L^2(0,T;H^1(\Omega)),\\
     \vr_\kappa\to\vr \quad \text{in\ } C([0,T]; L^{p}_{\rm{weak}}(\Omega)),\ 1\leq p<\infty,\\
     (1-\kappa)\vr_\kappa\to 0 \quad \text{strongly\ in\ }L^\infty(0,T;H^1(\Omega))\cap L^2(0,T; H^2(\Omega)), \\
      \vr_\kappa \rightharpoonup \vr  \hbox{ weakly$^*$ in } L^\infty(0,T;H^1(\Omega)), \qquad
      \vr_\kappa \rightharpoonup \vr   \hbox{ weakly  in } L^2 (0,T;H^2(\Omega)).
     \end{gathered}
     \end{equation}
and to pass to the limit in the continuity equation. Then the 
main difference with respect to limit process in the previous subsection concerns the strong convergence of the density. Note that the fourth convergence is useless because not uniform
with respect to $1-\kappa$.
  The strong convergence
$$\vr_\kappa\to \vr\quad \text{in\ } C([0,T]; L^2 (\Omega))$$
follows directly from the variant of the Aubin-Lions lemma, see f.i.  Corollary 4 in \cite{S87}
(using the uniform estimate of $\vr_\kappa$ in $L^\infty(0,T;H^1(\Omega)$ and the estimate
of $\partial_t\vr_\kappa \in L^2(0,T;H^{-1}(\Omega))$.
 Therefore we conclude on the strong convergence in $C([0,T];L^p(\Omega))$ for all
 $1\le p<+\infty$ using \eqref{conv2}${}_3$. Remark that we deduce 
 $$(\kappa \mu'(\vr_k))^{1/2} \nabla\vr_\kappa \to (\mu'(\vr))^{1/2} \nabla\vr 
     \hbox{ in } L^2(0,T;L^2(\Omega))$$
coming back the energy.      
    To pass to the limit in  the fifth the most demanding  and the only new term in comparison to \cite{BrEsSy07}  from the weak formulation of the momentum equation \eqref{weak_mom}. For this purpose one needs to check that $(1-\kappa)|\Grad\vr_\kappa|^2\to 0$ in $L^1((0,T)\times \Omega)$. But that is an immediate consequence of the fact that $\Grad\vr_\kappa$ is uniformly bounded in $L^2((0,T)\times\Omega)$. We therefore pass to the limit in \eqref{weak_mom} and find the
limit system
\eq{
&\intTO{\vr\vw\cdot\pt\vphi}+\intO{\vr(\vw-\kappa\Grad\vp(\vr))\otimes\vw:\Grad\vphi}\\
&\quad-2\intTO{\mu(\vr)A(\vw):\Grad\vphi} =-\intO{\vw^0\cdot\vphi(0)}\label{weak_mom1}}
holds for $\vphi\in (C^\infty((0,T)\times\Omega))^3$, s.t.  $\Div\vphi=0$ and $\vphi(T)=\vc{0}$
with $A(\vw)= (\nabla \vw - \nabla^t \vw)/2$. Using the regularity, we can write the weak
formulation in terms of $(\vr,\vu)$. This provides the convergence to a global weak solution
of the Kazhikhov-Smagulov-type system, for which the global solutions were proved in \cite{BrEsSy07}.

\section{Proof of Theorem \ref{T_main3} and an example.} \label{S:gen}
Below we present only the proof of the a priori estimate for more general viscosity and conductivity. The rest of the proof of existence of solutions would require minor modifications solely, so we skip this part.\\
Relation \eqref{gen_pm} leads to more general form of $\vw$, for the purposes of this section we denote
\eq{&\vw=\vu+ 2\Grad \tilde\vp(\vr),\\
&\Div\vw=0
\label{def_w2}}
and we define a new function $\tilde\mu(\vr)$ that satisfies
$$\vr\tilde\vp'(\vr)=\tilde\mu'(\vr).$$
Using this notation, system \eqref{Limref1} reads
\begin{equation}\label{mainor}
\begin{array}{c}
\pt\vr+ \Div(\vr \vu)  = 0,\\
\pt\lr{\vr \vw}+ \Div  (\vr \vu \otimes \vw) - 2 \Div (\mu(\vr) D(\vu)) + 2 \Div (\tilde\mu(\vr) \Grad^t \vu) + \Grad {\pi_1} =\vc{0},\\
\Div\vw = 0.
\end{array}
\end{equation}
We can again rewrite the momentum equation as
\eqh{\pt\lr{\vr \vw}+ \Div  (\vr \vu \otimes \vw) - 2 \Div ([\mu(\vr)-\tilde\mu(\vr)] D(\vu)) -2 \Div (\tilde\mu(\vr) A (\vu)) + \Grad {\pi_1} =\vc{0},}
and the energy estimate for this form is
\eq{\label{22}
\frac{1}{2}\Dt&\intO{\vr|\vw|^2}+ 2 \intO{[\mu(\vr)-\tilde\mu(\vr)] |D(\vu)|^2}\\
&+ 2 \intO{\tilde\mu(\vr)|A(\vu)|^2}+ 4 \intO{(\mu(\vr)-\tilde\mu(\vr))\Grad\vu:\Grad^2\tilde\vp(\vr)}=0.
}
Now, an extra estimate for $\Grad\tilde\vp(\vr)$ can be obtained by mimicking the steps leading to \eqref{b}, we have
\eq{\label{bb2}
\Dt\intO{\vr\frac{|\Grad\tilde\vp(\vr)|^2}{2}}-\intO{\tilde\mu(\vr)\Grad\vu:\Grad^2\tilde\vp(\vr)}-\intO{(\tilde\mu'(\vr)\vr-\tilde\mu(\vr))\Div\vu\lap\tilde\vp(\vr)}=0.
}
Next, we multiply \eqref{bb2} by the constant $4 \xi$ and we add it to \eqref{22} to find
\eq{\label{kaa}
\Dt\intO{\vr\lr{\frac{|\vw|^2}{2}+\xi\frac{|2\Grad\tilde\vp(\vr)|^2}{2}}}
+ 2 \intO{(\mu(\vr) -\tilde\mu(\vr))|D(\vu)|^2}+ 2\intO{\tilde\mu(\vr)|A(\vu)|^2}\\
+ 4 \intO{(\mu(\vr)-\tilde\mu(\vr)-\xi\tilde\mu(\vr))\Grad\vu:\Grad^2\tilde\vp(\vr)}
+ 2\xi\intO{(\tilde\mu'(\vr)\vr-\tilde\mu(\vr))|\Div\vu|^2}=0.
}
Therefore, using equivalence \eqref{Ddiv} we obtain
\eq{\label{kaa2}
\Dt&\intO{\vr\lr{\frac{|\vw|^2}{2}+\xi\frac{|2\Grad\tilde\vp(\vr)|^2}{2}}}
+2\intO{(\mu(\vr) -\tilde\mu(\vr)) |D(\vu)- \frac{1}{d} \Div\vu|^2}\\
&+2\intO{\tilde\mu(\vr)|A(\vu)|^2}
+4\intO{(\mu(\vr)-\tilde\mu(\vr)-\xi\tilde\mu(\vr))\Grad\vu:\Grad^2\tilde\vp(\vr)}\\
&+ 2 \intO{\lr{\xi(\tilde\mu'(\vr)\vr-\tilde\mu(\vr))+\frac{\mu(\vr) -\tilde\mu(\vr)}{d}}|\Div\vu|^2}=0.
}
Using the fact that $\intO{\Grad\vu:\Grad^2\tilde\vp(\vr)}=\intO{D(\vu):\Grad^2\tilde\vp(\vr)}$ we rewrite the fourth term
\eq{\label{integ0}
\Dt&\intO{\vr\lr{\frac{|\vw|^2}{2}+\xi\frac{|2\Grad\tilde\vp(\vr)|^2}{2}}}
+ 2\intO{(\mu(\vr) -\tilde\mu(\vr)) |D(\vu)- \frac{1}{d} \Div\vu|^2}\\
&+2\intO{\tilde\mu(\vr)|A(\vu)|^2}
+ 4\intO{(\mu(\vr)-\tilde\mu(\vr)-\xi\tilde\mu(\vr))\lr{D(\vu)-\frac{1}{d}\Div\vu \vc{I}}:\Grad^2\tilde\vp(\vr)}\\
&- 2\intO{\frac{\mu(\vr)-\tilde\mu(\vr)-\xi\tilde\mu(\vr)}{d}|\Div\vu|^2}\\
&+ 2 \intO{\lr{\xi(\tilde\mu'(\vr)\vr-\tilde\mu(\vr))+\frac{\mu(\vr) -\tilde\mu(\vr)}{d}}|\Div\vu|^2}=0.
}
and thus finally
\eq{\label{integ1}
\Dt&\intO{\vr\lr{\frac{|\vw|^2}{2}+\xi\frac{|2\Grad\tilde\vp|^2}{2}}}
+ 2\intO{\lr{\mu(\vr) -\tilde\mu(\vr)} |D(\vu)- \frac{1}{d} \Div\vu\,\vc{I}|^2}\\
&+ 2 \intO{\tilde\mu(\vr)|A(\vu)|^2}
+4\intO{(\mu(\vr)-\tilde\mu(\vr)-\xi\tilde\mu(\vr))\lr{D(\vu)-\frac{1}{d}\Div\vu\, \vc{I}}:\Grad^2\tilde\vp(\vr)}\\
&+2\intO{\xi\lr{\tilde\mu'(\vr)\vr+\frac{1-d}{d}\tilde\mu(\vr)}|\Div\vu|^2}=0.
}
Let us now denote
$$J_1(\vr)= \tilde\mu'(\vr)\vr+\frac{1-d}{d}\tilde\mu(\vr)
$$
$$J_2(\vr)= \mu(\vr) -\tilde \mu(\vr),
$$
From \eqref{c_gen} we know that there exists a positive constant $c$ such that
\eq{J_2(\vr) \ge c > 0, \qquad J_1(\vr)\geq c >0 \quad\hbox{ on } \quad[r,R]  \label{j1j2},}
therefore, the second and the last integral in \eqref{integ1} are non-negative, in fact we have
\eq{\label{integ11}
\Dt&\intO{\vr\lr{\frac{|\vw|^2}{2}+\xi\frac{|2 \Grad\tilde\vp(\vr)|^2}{2}}}
+ 2 \intO{J_2(\vr) |D(\vu)- \frac{1}{d} \Div\vu\, \vc{I}|^2}\\
&+2 \intO{\tilde\mu(\vr)|A(\vu)|^2}
+4 \intO{\lr{J_2(\vr)-\xi\tilde\mu(\vr)}\lr{D(\vu)-\frac{1}{d}\Div\vu\, \vc{I}}:\Grad^2\tilde\vp(\vr)}\\
&+2 \intO{\xi J_1(\vr)|\Div\vu|^2}=0.
}
So, in order to deduce uniform estimates from \eqref{integ1} we need to show that the penultimate term can be controlled by the positive contributions from the l.h.s. To this purpose let us
write
\eq{\nonumber
I= \intO{(J_2(\vr)-\xi\tilde\mu(\vr)) (D(\vu)- \frac{1}{d}\Div \vu\,  \vc{I}) :\Grad^2\tilde\vp(\vr)} \\
= 
\intO{\sqrt{J_2(\vr)}(D(\vu)-\frac{1}{d}{\Div \vu}\,\vc{I}):\frac{J_2(\vr)-\xi\tilde\mu(\vr)}{\sqrt{J_2(\vr)}}  \nabla\nabla \tilde\vp(\vr)
     },
}
then
$$I \le  \intO{\frac{(J_2(\vr)-\xi\tilde\mu(\vr))^2}{2 J_2(\vr)} |\Grad^2\tilde\varphi(\vr)|^2}+ \intO{\frac{J_2(\vr)}{2} |D(\vu)-\frac{1}{d}{\Div \vu}\, \vc{I}|^2 } 
$$
and the last term is absorbed by the l.h.s. of \eqref{integ11}.
Therefore, on account of equality
\eqh{\|\Grad  \Grad  \tilde\varphi(\vr)\|_{L^2(\Omega)} = \|\lap \tilde\varphi(\vr)\|_{L^2(\Omega)}}
it remains to assume that
\eqh{\max_{\vr\in[r,R]} \frac{(J_2(\vr)-\xi\tilde\mu(\vr))^2}{ 2J_2(\vr)}\leq \xi  \min_{\vr\in[r,R]} J_1(\vr)}
which is equivalent to \eqref{c_gen}. $\Box$

\medskip
Let us now consider a special case, satisfied for physical case of low Mach approximation for dense gases.
\begin{prop}\label{Example}
Assume that all the previous assumptions are satisfied and let \eqref{c_gen} be replaced by 
$$\mu(\vr)=\vr,\quad \tilde\mu(\vr)=\log\vr.$$
Then there exist a non-empty interval $[\tilde r,\tilde R]$ such that if
\eqh{0<\tilde r\leq \vr^0\leq \tilde R<0,}
there exists global in time weak solution to \eqref{Limref1}.
\end{prop}
\pf 
Conditions \eqref{c_gen} for $\vr=s=const.$ become
\begin{equation}
\begin{gathered}
c\leq s-\log s,\\
 \frac{(s-\log s-\xi\tilde\mu(s))^2}{2 \lr{s-\log s}}\leq \xi  c_1,\quad\text{where}\quad
 c_1\leq  1+\frac{1-d}{d}\log\tilde r.
\end{gathered}
 \label{c_s}
 \end{equation}
which means that $\xi\in\left[\xi^-,\xi^+\right]$, where
\eqh{
\xi^\pm=\frac{\lr{s-\log s}(\tilde\mu(s)+c_1)\pm \lr{s-\log s}\sqrt{c_1(2\tilde\mu(s)+c_1)}}{(\tilde\mu(s))^2}.
}
Now if $s$ is such that $\tilde\mu(s)>-\frac{c_1}{2}$, then at least $\xi^+$ is positive. In fact taking $\xi^0 =J_2(s)(\tilde\mu(s)+c_1)$, we always know that $\xi^0\in[\xi^-,\xi^+]$, so it satisfies conditions \eqref{c_s}.
Let us now note that the second condition in \eqref{c_s} is continuous with respect to $s$. Therefore from above considerations it follows that there exists a neighbourhood $(\tilde r,\tilde R)$ of $s$ in which \eqref{c_s} is satisfied. $\Box$\\

\smallskip
 
This example shows that Theorem \ref{T_main3}  generalizes the global in time existence of weak  solutions shown by {\sc P.--L. Lions} (see \cite{PLL}, Chapter 8.8) for the two-dimensional case to the three-dimensional case.

\section{Application to the model of the gaseous mixture}\label{S:mix}
As a particular application of the above theory, let us consider the system which describes the flow of a compressible $n$-component fluid in a domain $\Omega \subset \R^d$, $d=2,3$
\begin{equation}\label{1.1}
\begin{array}{c}
\pt\vr+\Div (\vr \vu) = 0,\\
\pt\lr{\vr\vu}+\Div (\vr \vu \otimes \vu) + \Div \vc{S} + \Grad p =\vc{0},\\
\pt\lr{\vr e}+\Div (\vr e\vu )+p\Div\vu +\Div\bf{Q}+ \vc{S}:\Grad\vu=0,\\
\pt\vr_k+\Div (\vr_k \vu)+ \Div (\vf_{k})  =  0,\quad k\in \{1,\ldots,n\},
\end{array}
\end{equation}
where the unknowns are  the species densities $\vr_k=\vr_k(t,x)$,  $k\in \{1,\ldots,n\}$, the total mass density $\vr=\vr(t,x)$, $\vr=\sumkN{\vr_{k}}$, the velocity vector field $\vu=\vu(t,x)$ and the absolute temperature $T=T(t,x)$. Further $\vc{S}$ denotes the viscous tensor, $p$ the internal pressure of the fluid,  $e$  the internal  energy, $\vc{Q}$ the heat flux, $m_k$ the molar mass of the $k$-th species and $\vf_{k}$ denotes the diffusion flux of the $k$-th species. 

\medskip

\noindent{\bf The equation of state.} For simplicity we consider a mixture of ideal gases with constant specific heats, i.e.
$$p=\sumkN p_k=\sumkN \frac{R\vr_k T}{m_k},\quad \vr e=\sumkN \vr_k e_k=\sumkN c_{vk}\vr_{k}T,$$
where $c_{vk}$ and $c_{pk}$ are the constant-volume and the constant-pressure specific heats related by
\eq{c_{pk}=c_{vk}+\frac{R}{m_k}\label{defcp}} 
and $m_k$ is the species $k$ molar mass.

\medskip

\noindent{\bf The diffusion fluxes.} We consider the so-called multicomponent diffusion without taking into account the Soret effect. The diffusion fluxes  may depend on $\vr, T,\vr_1,\ldots,\vr_n$ as follows
\begin{equation}
\vf_{k}=-C_0\sum_{l=1}^n C_{kl}\vc{d}_l,
\label{eq:diff}
\end{equation}
where $C_{0},\ C_{kl}$ are multicomponent flux diffusion coefficients;
by $\vd$ we denote the species $k$ diffusion force specified, in the absence of external forces, by the following relation
        \begin{equation}\label{eq:}
        \vd=\Grad\left({p_{k}\over p}\right)+\left({p_{k}\over p}-{\vr_{k}\over \vr}\right)\Grad\log{p}.
    \end{equation}
We assume that 
\begin{equation}\label{formD}
{C}_{kl}=D_{kl}\vr_{k},\quad k,l\in\{1,\ldots,n\},
\end{equation}
where $D_{kl}$ is symmetric and positive definite over the physical hyperplane ${U}^{\bot}$--the orthogonal complement of ${ U}=(1,\ldots,1)^{T}$, we refer the reader to \cite{VG}, Chapters 4, 7 for more details.

\medskip

\noindent{\bf The heat flux.} We neglect the transfer of energy due to species molecular diffusion---the so-called Dufour effect. The heat flux is thus in the following form
\begin{equation}\label{heat}
\vc{Q}=\sumkN c_{pk}T \vf_{k}-k\Grad T, \quad k>0,
\end{equation}
where the heat conductivity coefficient $k$ may depend smoothly on $\vr$ and $T$.
\medskip

\subsection{Low Mach system and simplifications}
When the Mach number is assumed to be small and for appropriate boundary conditions, system \eqref{1.1} has the following form
\begin{equation}\label{1.3}
\begin{array}{c}
\pt\vr+\Div (\vr \vu) = 0,\\
\pt\lr{\vr\vu}+\Div (\vr \vu \otimes \vu) + \Div \vc{S} + \Grad \Pi =\vc{0},\\
\sumkN c_{vk}\pt\lr{\vr_k T}+\sumkN c_{vk}\Div\lr{\vr_k T\vu}+P_0 \Div\vu=\Div\lr{k\Grad T}-{\sumkN\Div\lr{c_{pk}T\vf_k}}\\
\pt\vr_k+\Div (\vr_k \vu)+ \Div (\vf_{k})  = 0,\quad k\in \{1,\ldots,n\},\\
P_0=\frac{R\vr T}{\Ov{m}}
\end{array}
\end{equation}
where $P_0$ denotes the constant pressure, $\Ov{m}$ denotes the mean molar mass of the mixture given by
\eqh{
\frac{\vr}{\Ov{m}}=\sumkN\frac{\vr_k}{m_k}}
and the diffusion fluxes have the following simplified form
        \begin{equation}
        \vf_k=-C_0\sumlN C_{kl}\Grad\lr{\frac{p_k}{P_0 }}\label{eq:1}.
    \end{equation}

Below we propose some hypothesis that are adequate  for the Low Mach, constant Lewis number, binary mixture. 
\begin{enumerate}
\item[A1.] The number of species is $n=2$.
\item[A2.] The heat capacity at constant volume per unit mole of the mixture $C_{vk}$, are equal
\eqh{C_{vk}=m_k c_{vk}= C_v, \quad k=1,2.
\label{ass1}}
\item[A3.] The diffusion matrix coefficients are given by
\eqh{C=\frac{1}{\vr}\lr{\begin{array}{rr}
\vr_2&-\vr_1\\
-\vr_2& \vr_1
\end{array}},\quad C_0=c_0(\vr)\frac{m_1m_2}{\Ov{m}^2},
\label{Cform}}
where $c_0(\vr)>0$.
\item[A4.] The heat conductivity coefficient $k$ depends on the concentrations of the species and on $c_0$
\eqh{k=\frac{c_0(\vr)\lr{C_v+R}}{\Ov{m}}=\frac{c_0(\vr) C_p}{\Ov{m}},
\label{ass_symp_2}}
which is equivalent with assumption that the Lewis number is equal to 1.
\end{enumerate}

Under these simplifications system (\ref{1.3}) may be rewritten
\begin{equation*}\label{symp2}
\begin{array}{c}
\pt\vr+\Div (\vr \vu) = 0,\\
\pt\lr{\vr\vu}+\Div (\vr \vu \otimes \vu) + \Div \vc{S} + \Grad \Pi =\vc{0},\\
\displaystyle P_0 \Div\vu=\Div\lr{\frac{c_0 R}{\Ov{m}}\Grad T}-\sumkN\Div\frac{RT \vf_k}{m_k},\\
\pt\vr_k+\Div (\vr_k \vu)+ \Div \vf_{k}  = 0,\quad k\in \{1,2\},\\
\displaystyle P_0 =\frac{R\vr T}{\Ov{m}},
\end{array}
\end{equation*}
where 
\eqh{\vf_1=-c_0(\vr)\Grad Y_k,\quad \vf_2=-c_0(\vr)\Grad Y_2,\quad Y_i=\frac{\vr_i}{\vr},\ i=1,2.\label{fF}}
Let us now observe that
\eq{\Div\vu&=\Div\lr{\frac{c_0(\vr) R}{P_0 \Ov{m}}\Grad T-\sumkN\frac{RT \vf_k}{P_0 m_k}}\\
&=\Div\lr{\frac{c_0(\vr) R}{P_0 \Ov{m}}\Grad T+\frac{c_0(\vr)RT }{P_0 }\Grad\lr{\frac{1}{\Ov{m}}}}\\
&=\Div\lr{c_0(\vr)\Grad\lr{\frac{1}{\vr}}}.
\label{1.41}}
So, finally our system can be rewritten as
\begin{equation}\label{symp3}
\begin{array}{c}
\pt\vr+\Div (\vr \vu) = 0,\\
\pt\lr{\vr\vu}+\Div (\vr \vu \otimes \vu) - 2 \Div  (\mu(\vr) D(\vu)) - \Grad (\lambda(\vr)\Div\vu) + \Grad \Pi =\vc{0},\\
\Div\vu=\Div(c_0(\vr)\Grad\vr^{-1}),\\
\pt\vr_k+\Div (\vr_k \vu)-\Div(c_0(\vr)\Grad Y_k)  =  0,\quad k\in \{1,2\},\\
\displaystyle P_0=\frac{\vr T}{\Ov{m}}.
\end{array}
\end{equation}

Note that the three first equations of \eqref{symp2} form the low Mach number model obtained by {P.--L. Lions} in \cite{PLL}. The extension of this system to the combustion models was presented by {\sc A.~Majda} in \cite{Majda84} for binary mixture. In \cite{Embid87} {\sc P. Embid} proved the local-in-time existence of unique regular solution to zero Mach number equations for reacting multi-component compressible mixture.
Decoupling the first three equations from the subsystem of reaction-diffusion equations we obtain the Kazhikhov-Smagulov type system. The first global existence result for particular choice of viscosity coefficient
$$\mu(\vr)=\frac{c_0}{2}\log\vr$$
and for the initial density bounded away from zero is due to {\sc D. Bresch} {\it et al.} \cite{BrEsSy07}. 
Concerning the reaction-diffusion equations, the existence analysis for such system in case of diffusion (\ref{eq:diff}-\ref{eq:}) and with the diffusion matrix as in \eqref{Cform}  was recently performed by {\sc P. B. Mucha}, {\sc M. Pokorn\'y} and {\sc E. Zatorska }in \cite{MPZ}. Their result holds  under certain regularity assumptions imposed on $\vr$ and $\vu$. 
	
\bigskip

\subsection{Existence of solutions}
Our strategy is to first resolve  the combustion system looking for $\vr,\vv$ and then use them to handle the reaction-diffusion equations.
For simplicity, we assume that $\Omega=\mathbb{T}^3$ and we supplement system \eqref{symp3} by the initial conditions
\begin{equation}
\begin{array}{c}
\vr(0,x)=\vr^0(x),\quad \vu(0,x)=\vu^0(x)\quad x\in\Omega,\\
\vr_k(0,x)=\vr_k^0(x),\quad0< r_k\leq \vr_k^0(x)\leq R_k<\infty,\quad k=1,2,\\
 \vr_1^0+\vr_2^0=\vr^0.
\label{init1}
\end{array}
\end{equation}

\begin{rmk}
We will restrict to the case when the initial density $\vr^0$ is bounded away from vacuum. This is due to nonlinearity in the continuity equation. For example  for a model of pollutant (studied in  {\rm \cite{BrEsSy07}}) there exists also existence result for $\vr^0\geq0$ {\rm \cite{Sy2005}}.
\end{rmk}
Taking 
$$\kappa=c_0(\vr)/\tilde c_0(\vr),\quad\text{and}\quad \vp'(\vr)=\tilde c_0(\vr)\vr^{-2}$$ 
we obtain from relation \eqref{rel1} that 
$$\mu'(\vr)=\tilde c_0(\vr)\vr^{-1}$$ 
and the three first equations of \eqref{symp3} give exactly system \eqref{main1} studied at the beginning of this paper. 
\begin{rmk}
For example for $c_0(\vr)=\kappa$, condition \eqref{important} 
is satisfied provided the bound from above for the density $R$ is sufficiently small, i.e.
$$1-\lr{1-\frac{1}{d}}\log R\geq0. $$
\end{rmk}
In general we have to assume that
   \eqh{
   \min_{\vr\in[r,R]} \lr{\frac{1-d}{d}\int_r^\vr \tilde c_0(s)s^{-1}{\rm d}s+\tilde c_0(\vr)}\geq c > 0.
   }
With this information at hand we may combine the elements of the proof from  \cite{MPZ1} in order to construct weak solution to the system of reaction-diffusion equations for the species. The difference is that now we consider a general form of diffusion matrix $C$ but, on the other hand the form of diffusion flux $\vf_k$  is reduced to \eqref{eq:1}.
Thanks to this, the species equations reads
\eqh{\pt(\vr_\ep Y_k)+\Div (\vr_\ep Y_k \vu_\ep)-\Div(c_0(\vr_\ep)\Grad Y_k) = 0,\quad k=1,2,}
where $Y_k=\frac{\vr_k}{\vr}$, $\vr_\ep,\ \vu_\ep$ denotes the convolution of $\vr, \vu$ with a standard regularizing kernel. 
Due to maximum principle for $\vr$ \eqref{max_vr}, this is system with semi-linear parabolic equations which are weakly coupled. The existence of  strong solutions and the limit passage $\ep\to0$ is then straightforward. 

\section{Application to the ghost effect system}\label{S:ghost}
  In this section, we give some comments regarding a particular case of
 ghost system that we can find for instance in a recent paper of 
 {\sc C.D. Levermore, W. Sun,  K. Trivisa} \cite{LeSuTr}
  \begin{equation*}
\begin{array}{c}
\vr T=1, \qquad 
\partial_t \vr + \Div (\vr \vu) = 0, \\
\pt (\vr \vu) + \Div (\vr \vu \otimes \vu)
         +\nabla P^*= - \Div\vc{\Sigma}-\Div\tilde{\vc{\Sigma}} ,\\
\frac{5}{2}\Div\vu =\Div\lr{k(T)\Grad T},
\end{array}
\end{equation*} 
where $k(T)$ is the heat conductivity coefficient $k>0$, while the two parts of stress tensor $\vc{\Sigma}$ and $\tilde{\vc{\Sigma}}$ are defined as follows
\eqh{
\vc{\Sigma}&=\mu(T)\lr{\Grad\vu+\Grad^t\vu-\frac{2}{3}\Div\vu \vc{I}},\\
\tilde{\vc{\Sigma}}&=\tau_1(\vr,T)\lr{\Grad^2T-\frac{1}{3}\lap T\vc{I}}+\tau_2(\vr,T)\lr{\Grad T\otimes\Grad T-\frac{1}{3}|\Grad T|^2\vc{I}},
}
where $\tau_1,\tau_2$ are transport coefficients with $\tau_1>0$.
Now, let us consider a particular form of heat-conductivity coefficient, such that
$$\frac{2}{5}\frac{k\lr{{\vr}^{-1}}\Grad\vr}{\vr^2}=\kappa\Grad\log\vr$$
and let us choose $\tau_1=c\vr^2$, $\tau_2=-c\vr^3$ , then we obtain
\eqh{
\tau_1\Grad^2 T&= c\vr^2\Grad^2\lr{\frac{1}{\vr}}=-c\vr\Grad^2\log\vr+c\vr^{-1}\Grad\vr\otimes\Grad\vr\\
\tau_1\lap T&=-c\vr\lap\log\vr+c\frac{|\Grad\vr|^2}{\vr}
}
and therefore
\eqh{\tau_1(\vr, T)\lr{\Grad^2 T-\frac{1}{3}\lap T\vc{I}}&=-c\vr\Grad^2\log\vr
+c\vr^{-1}\Grad\vr\otimes\Grad\vr-\frac{c}{3}\frac{|\Grad\vr|^2}{\vr}\vc{I}+\frac{c}{3}\vr\lap\log\vr\vc{I},\\
\tau_2(\vr, T)\lr{\Grad T\otimes\Grad T-\frac{1}{3}|\Grad T|^2\vc{I}}&=-c\vr^{-1}\Grad\vr\otimes\Grad\vr+\frac{c}{3}\frac{|\Grad\vr|^2}{\vr}\vc{I}
}
and hence
\eqh{\tilde{\vc{\Sigma}}=-c\vr\Grad^2\log\vr+\frac{c}{\kappa}\vr\Div\vu\vc{I}}
and the last term, after taking the divergence, may be absorbed by a part of $\vc{\Sigma}$ or put to the pressure.
Finally, the ghost system can be rewritten as follows
 \begin{equation}\label{Ghost}
\begin{array}{c}
\partial_t \vr + \Div (\vr \vu) = 0. \\
\pt (\vr \vu) + \Div (\vr \vu \otimes \vu)
         +\nabla \pi -  2 \Div (\mu(\vr) D(\vu))   - c\, \Div(\vr \nabla\nabla \log\vr) = 0,\\
\Div \vu = - 2 \kappa \Delta \log\vr.
\end{array}
\end{equation} 
Let us focus on the specific form of viscosity coefficient $\mu(\vr)= \bar{\mu}\vr$,
 thanks to the Bohm identity we have
$$\Div(\vr\nabla\nabla \log\vr) = 2\vr \nabla (\frac{1}{\sqrt\vr}\Delta \sqrt\vr).$$
   System \eqref{Ghost} with above relation corresponds to  the low Mach version with large heat-release of the quantum Navier-Stokes  system studied by A. {\sc J\"ungel} in \cite {Ju},  by M. {\sc Gisclon} and I. {\sc Violet} in \cite{GiVi}, see also work by B. Haspot \cite{Ha14}.

First we check that  system \eqref{Ghost} has the following energy  estimate
\eq{\label{estimghost}
\Dt&\intO{\vr\lr{\frac{|\vw|^2}{2}   +[(1-\kappa)\kappa+c] \frac{|2\Grad\log\vr|^2}{2}}} 
 + 2(1-\kappa)\bar{\mu}\intO{\vr |D( \vu)|^2}  \\
&+ 2\kappa\bar{\mu}\intO{\vr |A(\vu)|^2} + c\kappa \bar{\mu} \intO{\vr|\Grad^2\log\vr|^2}=0
}
with $\vw=\vu + 2 \kappa \bar{\mu} \nabla\log\vr$.
Moreover we have from the maximum principle
$$r\le \vr \le R$$
and the above estimate implies that for  $0<\kappa< 1$ we have $\vr \in L^2(0,T;H^2(\Omega)).$
This is the final argument needed to justify global in time existence of weak solutions to the ghost effect system.  System \eqref{Ghost}, in comparison to the system studied in the first part of the paper, has a new term which needs to be handled, namely 
\eqh{\Div(\vr \nabla\nabla \log\vr)=\Grad\lap\vr-\Div\lr{\Grad\sqrt{\vr}\otimes\Grad\sqrt{\vr}}.}
The difficulty to pass to the limit is the last component  of this expressionn which asks for strong convergence of $\Grad\sqrt{\vr}$ in $L^2((0,T)\times\Omega)$.
This is given using the regularity obtained by \eqref{estimghost} and the mass equation
using the standard  Aubin-Lions lemma.

\medskip

\begin{rmk}
Similar regularity argument was used by {\sc I. Kostin, M. Marion, R.~Texier-Picard} and {\sc V.A. Volpert}\cite{KMT-PV} to prove the global existence of solutions to the incompressible Korteweg system but with an additional diffusion in the mass (concentration) equation.
\end{rmk}
\begin{rmk}
Estimate \eqref{estimghost} is not really helpful to perform the limit $\kappa\to 0$. The global existence of weak solutions to the system with $\kappa=0$ is in fact an open problem due to the lack of convergence in the nonlinear capillary term, namely in the quantity $\nabla\sqrt\vr \otimes\nabla\sqrt\vr$.
\end{rmk}

\bigskip

\noindent {\bf Acknowledgments.}
 The authors thank the referee for his/her valuable comments which improves the quality of  the paper. The first and the third author acknowledge the support from the ANR-13-BS01-0003-01 project DYFICOLTI.   
 The second and the third author acknowledge the Post-Doctoral support of Ecole Polytechnique.
 The third author was also supported by MN grant  IdPlus2011/000661 and by the fellowship
START of the Foundation for Polish Science.

\end{document}